\documentclass[12pt]{article}
\usepackage{amssymb}
\usepackage{amsfonts}
\usepackage{amsmath}

\input{tcilatex}

\begin{document}

\title{ACTION TYPE GEOMETRICAL EQUIVALENCE OF REPRESENTATIONS OF GROUPS.}
\author{}
\maketitle

\begin{center}
{\LARGE B. Plotkin}$^{\sharp }${\LARGE , A.Tsurkov}$^{\flat }$\bigskip

$^{\sharp }$\textit{\ Institute of Mathematics}

\textit{Hebrew University, Givat Ram, 91904}

\textit{Jerusalem, 91904, Israel}

borisov@math.huji.ac.il

\bigskip

and

\bigskip

$^{\flat }$\textit{\ Department of Mathematics and Statistics}

\textit{Bar Ilan University}

\textit{Ramat Gan, 52900, Israel.}

tsurkoa@macs.biu.ac.il

\bigskip
\end{center}

\begin{abstract}
In the paper we prove (\textbf{Theorem 8.1})  that \textit{there
exists a
continuum of non isomorphic simple modules over }$KF_{2}$\textit{, where }$%
F_{2}$\textit{\ is a free group with }$2$\textit{\ generators}
(compare with [Ca] where a continuum of non isomorphic simple
$2$-generated groups is constructed).  Using this fact
 we give an example of a non action type logically
Noetherian representation (\textbf{Section} 9).

In general, the topic of this paper is the action type algebraic
geometry of representations of groups. For every variety of algebras
$\Theta $ and every algebra $H\in \Theta $ we can consider an
algebraic geometry in $\Theta $ over $H$. Algebras in $\Theta $ may
be many sorted (not necessarily one sorted) algebras. A set of sorts
$\Gamma $ is fixed for each $\Theta $. This theory can be applied to
the variety of representations of groups over fixed commutative ring
$K$ with unit. We consider a representation as two sorted algebra
$\left( V,G\right) $, where $V$ is a $K$-module, and $G$ is a group
acting on $V$.
In the action type algebraic geometry of representations of groups algebraic
sets are defined by systems of action type equations and equations in the
acting group are not considered. This is the special case, which cannot be
deduced from the general theory (see \textbf{Corollary} from \textbf{%
Proposition 3.5}, \textbf{Corollary 2 }from \textbf{Proposition 4.2} and
\textbf{Remark 5.1}). In this paper the following basic notions are
introduced: action type geometrical equivalence of two representations,
action type quasi-identity in representations, action type quasi-variety of
representations,\textbf{\ }action type Noetherian variety of
representations, action type geometrically Noetherian representation, action
type logically Noetherian representation. \textbf{Proposition 6.2}, and
\textbf{Corollary} from \textbf{Proposition 6.3} provide examples of action
type Noetherian variety of representations and action type geometrically
Noetherian representations. In \textbf{Corollary 2} from \textbf{Theorem 5.1}
the approximation-like criterion for two representations to be action type
geometrically equivalent is proved. This criterion is similar to the
approximation criterion for two algebras to be geometrically equivalent in
regular sense ([PPT]). \textbf{Theorem 6.2} gives a criterion for a
representation to be action type logically Noetherian. This criterion is
formulated in terms of an action type quasi-variety generated by a
representation (compare with [Pl4]). In \textbf{Corollary 2} from \textbf{%
Theorem 7.1} we consider a Birkhoff-like description [Bi] of an action type
quasi-variety generated by a class of representations. An example of a non action type
logically Noetherian representation allows to build an ultrapower of a
non action type logically Noetherian representation, which has the same
action type quasi-identities but is not action type geometrically equivalent
to the original representation ((\textbf{Corollary} from \textbf{Theorem 9.1}%
). This result is parallel to the corresponding theorem for groups [MR].
\end{abstract}

\centerline{\bf{ Introduction.}}

\bigskip

In this paper we consider the action type algebraic geometry of
representations of groups. General references for universal
algebraic geometry, i.e. the geometry associated with varieties of
algebras  are [BMR], [MR], [Pl1--Pl4]. First notions in the
algebraic geometry of representations of groups were defined in
[Pl5]. We outline them in the introduction and consider in detail in
the sequel.

We consider only right side modules and throughout the paper "module" means
a "right side module". Let $K$ be a commutative ring with unit, $G$ be a
group, $V$ be a $K$-module, and $KG$ be the group ring over the group $G$. $%
\left( V,G\right) $ is a representation of the group $G$ if $V$ is a $KG$%
-module. This is equivalent to the existence of the group homomorphism $\rho
:G\rightarrow \mathrm{Aut}_{K}\left( V\right) $ and the ring homomorphism $%
\rho :KG\rightarrow \mathrm{End}_{K}\left( V\right) $ (in the paper the
homomorphism of groups $\varphi :G\rightarrow H$, and the corresponding ring
homomorphism $\varphi :KG\rightarrow R,$ where $R\supseteq KH$ are denoted
by the same letter). The multiplication of elements of the module $V$ by
elements of $G$ and $KG$ is denoted by $\circ $ and other similar symbols,
and is called the action of the group $G$ (ring $KG$) on elements of the
module $V$. The variety of representations of groups over the fixed
commutative ring $K$ we denote $Rep-K$ ([PV]).

The homomorphism of two representations $\left( \alpha ,\beta \right)
:\left( V,G\right) \rightarrow \left( W,H\right) $, is the pair $\left(
\alpha ,\beta \right) $ of two homomorphisms where $\alpha :V\rightarrow W$
is the homomorphism of $K$-modules and $\beta :G\rightarrow H$ is the
homomorphism of groups subject to condition $\left( v\circ g\right) ^{\alpha
}=v^{\alpha }\circ g^{\beta }$, for every $v\in V$ and $g\in G$. If $\left(
V,G\right) $ is a representation, $V_{0}\leq V$ is a $K$-submodule of $V$, $%
G_{0}\leq G$ is a subgroup of $G$ and $V_{0}$ is a $KG_{0}$-submodule, then
we say that $\left( V_{0},G_{0}\right) $ is a subrepresentation of $\left(
V,G\right) $ (denoted by $\left( V_{0},G_{0}\right) \leq \left( V,G\right) $%
).

If $\left( \alpha ,\beta \right) :\left( V,G\right) \rightarrow \left(
W,H\right) $ is a homomorphism of representations, $V_{0}=\ker \alpha $, $%
G_{0}=\ker \beta $, then we denote $\ker \left( \alpha ,\beta \right)
=\left( V_{0},G_{0}\right) $. We have $\ker \left( \alpha ,\beta \right)
\leq \left( V,G\right) $, $G_{0}\trianglelefteq G$, $V_{0}$ is a $KG$-module
and $G_{0}$ acts trivially on the $V/V_{0}$. On the other hand, if $\left(
V_{0},G_{0}\right) \leq \left( V,G\right) $ is a subrepresentation, which
satisfies the conditions 1) $G_{0}\trianglelefteq G$, 2) $V_{0}$ is a $KG$%
-module and 3) $G_{0}$ acts trivially on the $V/V_{0}$, then  one
can define the action of $G/G_{0}$ on the $V/V_{0}$ by the rule:
$\left( v+V_{0}\right) \circ \left( gG_{0}\right) =v\circ g+V_{0}$
($v\in V$, $g\in G $). Then $\left( V/V_{0},G/G_{0}\right) $ is the
representation and the pair of natural homomorphisms $\alpha
:V\rightarrow V/V_{0}$, $\beta :G\rightarrow G/G_{0}$ is the
homomorphism of representations $\left( \alpha
,\beta \right) :\left( V,G\right) \rightarrow \left( V/V_{0},G/G_{0}\right) $%
. Subrepresentation $\left( V_{0},G_{0}\right) \leq \left( V,G\right) $,
which satisfies the conditions 1), 2) and 3) is called a normal
subrepresentation (denoted by $\left( V_{0},G_{0}\right) \trianglelefteq
\left( V,G\right) $).

Free objects in the $Rep-K$ are representations $\left( XKF\left( Y\right)
,F\left( Y\right) \right) $, where $F\left( Y\right) $ is the free group
with the set of free generators $Y$, $KF\left( Y\right) $ the group ring
over this group, and $XKF\left( Y\right) =\bigoplus\limits_{x\in X}xKF\left(
Y\right) $ is the free $KF\left( Y\right) $-module with the basis $X$. We
denote the representation $\left( XKF\left( Y\right) ,F\left( Y\right)
\right) $ as $W\left( X,Y\right) $. Below in this paper we suppose (if we do
not say anything else specifically) that $X$ and $Y$ are finite subsets of
the countable sets $X_{0}$ and $Y_{0}$ respectively.

In the variety $Rep-K$ we can consider subvarieties $\Theta $
defined simultaneously by a set of identities in acting groups
$\left\{ f=1\mid f\in F_{\infty }\right\} =H$, where $F_{\infty }$
is the free group of the countable rank and by a set of identities
of the form $\left\{ w=0\mid w\in
\mathbb{N}
KF_{\infty }\right\} =A$, where $%
\mathbb{N}
KF_{\infty }$ is the free $KF_{\infty }$-module of the countable rank.
Elements of $A$ are identities which describe action of groups on modules.
These identities are called action type identities. In other words $\Theta $
is the set of representations $(V,G)\in Rep-K$ such that $\forall \left(
\alpha ,\beta \right) \in \mathrm{Hom}\left( \left(
\mathbb{N}
KF_{\infty },F_{\infty }\right) ,\left( V,G\right) \right) $ holds $%
(f^{\beta }=1)\wedge (w^{\alpha }=0)\quad \forall f\in H,\ \forall
w\in A.$
Subvarieties of this kind are called in [PV] "bivarieties". For every $%
\left( V,G\right) \in \Theta $ we can consider the set of group identities
satisfied by $(V,G)$%
\begin{equation*}
Id_{gr}\left( V,G\right) =\left\{ f\in F_{\infty }\mid \forall \left( \alpha
,\beta \right) \in \mathrm{Hom}\left( \left(
\mathbb{N}
KF_{\infty },F_{\infty }\right) ,\left( V,G\right) \right) \left( f^{\beta
}=1\right) \right\} ,
\end{equation*}%
and the set of action-type identities satisfied by $(V,G)$
\begin{equation*}
Id_{a.t.}\left( V,G\right) =\left\{ w\in
\mathbb{N}
KF_{\infty }\mid \forall \left( \alpha ,\beta \right) \in \mathrm{Hom}\left(
\left(
\mathbb{N}
KF_{\infty },F_{\infty }\right) ,\left( V,G\right) \right) \left( w^{\alpha
}=0\right) \right\} .
\end{equation*}%
$Id_{gr}\left( V,G\right) \supseteq H$ and $Id_{a.t.}\left( V,G\right)
\supseteq A$ holds. Clearly, $\left( Id_{a.t.}\left( V,G\right)
,Id_{gr}\left( V,G\right) \right) $ is the normal subrepresentation of the
representation $\left(
\mathbb{N}
KF_{\infty },F_{\infty }\right) $. Denote%
\begin{equation*}
Id_{gr}\Theta =\dbigcap\limits_{\left( V,G\right) \in \Theta }Id_{gr}\left(
V,G\right) ,Id_{a.t.}\Theta =\dbigcap\limits_{\left( V,G\right) \in \Theta
}Id_{a.t.}\left( V,G\right) .
\end{equation*}%
$\left( Id_{a.t.}\Theta ,Id_{gr}\Theta \right) $ is also the normal
subrepresentation of the $\left(
\mathbb{N}
KF_{\infty },F_{\infty }\right) $ and $\dbigcap\limits_{\left( V,G\right)
\in \Theta }Id_{gr}\left( V,G\right) \supseteq H$, $Id_{a.t.}\Theta
\supseteq A$. By [PV], action type identities (elements of $Id_{a.t.}\Theta $%
) can be reduced to the identities in the cyclic module $\left\{ x\right\}
KF_{\infty }\cong KF_{\infty }$, i.e., to the identities of the form $x\circ
u\left( y_{1},\ldots ,y_{n}\right) =0,$ where $y_{1},\ldots ,y_{n}$ are some
generators of $F_{\infty }$, and $u\in KF_{\infty }$.

Example. The identity
\begin{equation*}
x\circ (y_{1}-1)(y_{2}-1)\ldots (y_{n}-1)
\end{equation*}%
defines the $n$-stable variety of group representations. This variety we
denote by $\mathfrak{S}^{n}$.

Denote%
\begin{equation*}
Id_{a.t.}\Theta \cap XKF\left( Y\right) =Id_{a.t.}\left( \Theta ,X,Y\right)
,Id_{gr}\Theta \cap F\left( Y\right) =Id_{gr}\left( \Theta ,Y\right) ,
\end{equation*}%
\begin{equation*}
Id\left( \Theta ,X,Y\right) =\left( Id_{a.t.}\left( \Theta ,X,Y\right)
,Id_{gr}\left( \Theta ,Y\right) \right)
\end{equation*}%
($Id\left( \Theta ,X,Y\right) $ is a normal subrepresentation of $W(X,Y)$)
\begin{equation*}
XKF\left( Y\right) /Id_{a.t.}\left( \Theta ,X,Y\right) =E_{\Theta }(X,Y),\
F\left( Y\right) /Id_{gr}\left( \Theta ,Y\right) =F_{\Theta }\left( Y\right)
.
\end{equation*}%
Then%
\begin{equation*}
W_{\Theta }(X,Y)=W(X,Y)/Id\left( \Theta ,X,Y\right) =\left( E_{\Theta
}(X,Y),F_{\Theta }\left( Y\right) \right)
\end{equation*}%
is the free representation in the variety of representations $\Theta $
(relatively free representation). Below in this paper we suppose (if we do
not say anything else specifically) that $\Theta $ is a subvariety of $Rep-K$%
.

Let $\left( V,G\right) $ be a fixed representation, such that $\left(
V,G\right) \in \Theta $. We consider \textit{affine spaces} of finite rank
over the $\left( V,G\right) $ in the variety $\Theta $. These are the sets $%
\mathrm{Hom}\left( W_{\Theta }\left( X,Y\right) ,\left( V,G\right) \right) $%
, where $W_{\Theta }\left( X,Y\right) $ is a finitely generated free
representation. We assume that the sets $X$ and $Y$ are finite
besides of the few explicitly pointed places.

We have two kinds of equations in the algebraic geometry over
representations: equations in the acting group of the form $f=1$, where $%
f\in F_{\Theta }\left( Y\right) $, and the action type equations of the form
$w=0$, where $w\in E_{\Theta }(X,Y)$. Action type equations describe action
of the group on a module.

In action type algebraic geometry of representations we consider only action
type equations. If $T\subset E_{\Theta }(X,Y)$ is a set of these equations,
it defines  the algebraic set%
\begin{equation*}
T_{\left( V,G\right) }^{\nabla }=\left\{ \left( \alpha ,\beta \right) \in
\mathrm{Hom}\left( W_{\Theta }\left( X,Y\right) ,\left( V,G\right) \right)
\mid \ker \alpha \supset T\right\} .
\end{equation*}%
over a representation $\left( V,G\right) $. If $A\subset
\mathrm{Hom}\left( W_{\Theta }\left( X,Y\right) ,\left( V,G\right)
\right) $ is a set of points of the affine space, then we have the
"ideal" of action type equations (in fact the $KF_{\Theta }\left(
Y\right) $-submodule in $E_{\Theta }(X,Y)$), defined by the set $A$:%
\begin{equation*}
A_{\left( V,G\right) }^{\nabla }=\bigcap\limits_{\left( \alpha ,\beta
\right) \in A}\ker \alpha .
\end{equation*}%
Now we can consider the action type $\left( V,G\right) $-closure of a set of
action type equations $T\subset E_{\Theta }(X,Y)$:%
\begin{equation*}
T_{\left( V,G\right) }^{\nabla \nabla }=\bigcap\limits_{\left( \alpha ,\beta
\right) \in T_{\left( V,G\right) }^{\nabla }}\ker \alpha
\end{equation*}%
and the action type $\left( V,G\right) $-closure of a set of points $%
A\subset \mathrm{Hom}\left( W_{\Theta }\left( X,Y\right) ,\left( V,G\right)
\right) $:
\begin{equation*}
A_{\left( V,G\right) }^{\nabla \nabla }=\left\{ \left( \alpha ,\beta \right)
\in \mathrm{Hom}\left( W_{\Theta }\left( X,Y\right) ,\left( V,G\right)
\right) \mid \ker \alpha \supset A_{\left( V,G\right) }^{\nabla }\right\} .
\end{equation*}

\textbf{Definition 0.1.} \textit{Representations }$\left( V_{1},G_{1}\right)
,\left( V_{2},G_{2}\right) \in \Theta $\textit{\ are called }\textbf{action
type}\textit{\ }\textbf{geometrically equivalent}\textit{\ (}denoted\textit{%
\ }$\left( V_{1},G_{1}\right) \sim _{a.t.}\left( V_{2},G_{2}\right) $\textit{%
) if }$T_{\left( V_{1},G_{1}\right) }^{\nabla \nabla }=T_{\left(
V_{2},G_{2}\right) }^{\nabla \nabla }$\textit{\ for every }$X$\textit{\ and }%
$Y$\textit{\ and for every set }$T\subset E_{\Theta }(X,Y)$\textit{.}

By \textbf{Proposition 4.3} this definition is correct, i.e., action type\
geometric equivalence of representations does not depend on  subvariety $%
\Theta $.

\textbf{Definition 0.2.}\textit{\ The universal logic formula of the form}%
\begin{equation}
\left( \bigwedge\limits_{i=1}^{n}\left( w_{i}=0\right) \right) \Rightarrow
\left( w_{0}=0\right) ,  \tag{0.1}
\end{equation}%
\textit{where }$\left\{ w_{0},w_{1},\ldots ,w_{n}\right\} \subset XKF\left(
Y\right) $\textit{, is called an }\textbf{action type quasi-identity}\textit{%
.}

We say that a representation $\left( V,G\right) $ is fulfilled on the
formula (0.1) and denote:%
\begin{equation*}
\left( V,G\right) \vDash \left( \left( \bigwedge\limits_{i=1}^{n}\left(
w_{i}=0\right) \right) \Rightarrow \left( w_{0}=0\right) \right)
\end{equation*}%
if $w_{0}^{\alpha }=0$ for every $\left( \alpha ,\beta \right) \in \mathrm{%
Hom}\left( W\left( X,Y\right) ,\left( V,G\right) \right) $, such that $%
w_{i}^{\alpha }=0$ for every $i\in \left\{ 1,\ldots ,n\right\} $.

Also we can consider the \textbf{infinite action type "quasi-identity"}:%
\begin{equation}
\left( \bigwedge\limits_{i\in I}\left( w_{i}=0\right) \right) \Rightarrow
\left( w_{0}=0\right) ,  \tag{0.1'}
\end{equation}%
where $\left\{ w_{0}\right\} \cup \left\{ w_{i}\mid i\in I\right\} \subset
XKF\left( Y\right) $, $\left\vert I\right\vert \geq \aleph _{0}$. This is
not a logic formula in the usual sense, but we can interpret its meaning by
the rule: a representation $\left( V,G\right) $ satisfies (0.1') if $%
w_{0}^{\alpha }=0$ for every $\left( \alpha ,\beta \right) \in \mathrm{Hom}%
\left( W\left( X,Y\right) ,\left( V,G\right) \right) $, such that $%
w_{i}^{\alpha }=0$ for every $i\in I$.

\textbf{Definitions 0.3. }\textit{A representation }$\left( V,G\right) \in
\Theta $ \textit{is called }\textbf{action type geometrically Noetherian}%
\textit{\ if for every sets }$X$\textit{\ and }$Y$\textit{\ and every set }$%
T\subset E_{\Theta }(X,Y)$\textit{, there is a finite set }$T_{0}\subset T$%
\textit{, such that }$\left( T_{0}\right) _{\left( V,G\right) }^{\nabla
}=T_{\left( V,G\right) }^{\nabla }$\textit{.}

\textit{A representation }$\left( V,G\right) \in \Theta $\textit{\ is called
}\textbf{action type logically Noetherian}\textit{\ if for every sets }$X$%
\textit{\ and }$Y$\textit{,\ every set }$T\subset E_{\Theta }(X,Y)$\textit{\
and every }$w\in T_{\left( V,G\right) }^{\nabla \nabla }$\textit{, there is
a finite set }$T_{0}\subset T$\textit{, such that }$w\in \left( T_{0}\right)
_{\left( V,G\right) }^{\nabla \nabla }$\textit{.}

Also by \textbf{Proposition 4.3}, action type geometric
Noetherianity and action type logic Noetherianity of representation
does not depend on subvariety $\Theta $.

The paper is organized as follows. We start with two auxiliary sections. For
the sake of completeness we recall in \textbf{Section 1} some basic
definitions and constructions for modules and representations of groups
which will be needed later. In \textbf{Section 2} we consider operators on
classes of algebras. Some of these operators act specifically on classes of
representations of groups and can be found in [PV]. In this paper we
continue to study the properties of these operators.

In \textbf{Section 3} we study the basic notions related to algebraic
geometry in representations of groups. We distinguish two kinds of
equations: equations in the acting group and action type equations. The main
concepts of the action type algebraic geometry of representations, i.e. the
geometry determined by action type equations, are defined in \textbf{Section
4}.\textbf{\ Section 5} deals with the notion of action type geometrical
equivalence of representations. In \textbf{Corollary 2} from \textbf{Theorem
5.1} an approximation-like criterion for two representations of groups to be
action type geometrically equivalent is presented.

The notions of Noetherian variety of algebras and geometrically
(logically) Noetherian algebra play an important role in the theory
(see, [Pl3--Pl4]). The corresponding notions of action type
Noetherian variety of representations and  action type geometrically
(logically) Noetherian representation are discussed in
\textbf{Section 6}. \textbf{Theorems 6.1} and \textbf{6.2} establish
relations between geometrical and logical properties of logically
Noetherian representations. Two examples are presented: the
$n$-stable variety of representation over the Noetherian ring $K$\
is the action type Noetherian variety for every $n\in
\mathbb{N}
$ (\textbf{Proposition 6.2}), and every finite dimension representation over
the field $K$\ is action type geometrically Noetherian (\textbf{Corollary}
from the \textbf{Proposition 6.3}).

R.Gobel and S. Shelah ([GSh]) proved that there is a non logically
Noe\-therian group. A.Myasnikov and V.Remeslennikov [MR] proved that
for every non logically Noetherian group there exists an ultrapower
of this group, which, of course, has the same quasi-identities as
the original group, but is not geometrically equivalent to the
original group. Our target in the three final sections is to prove a
similar result in the action type algebraic geometry of
representations. In \textbf{Section 7}\ we consider action type
quasi-varieties of representations, i.e. quasi-varieties of
representations, defined by action type quasi-identities. We give a
description of the action type quasi-variety generated by a class of
representations in terms of operators on classes of representations.

In \textbf{Section 8} we prove \textbf{Theorem 8.1}: \textit{There exists a
continuum of non isomorphic simple modules over }$KF_{2}$\textit{, where }$%
F_{2}$\textit{\ is a free group with }$2$\textit{\ generators} ($K$ is a
countable field). This theorem is similar to the result of R.Camm [Ca]:
there is a continuum of non isomorphic simple $2$-generated groups. The
latter theorem has been used by R.Gobel and S. Shelah in the construction of
a non logically Noetherian group. We use \textbf{Theorem 8.1} in the \textbf{%
Section 9} (\textbf{Theorem 9.1}) for the construction of a non action type
logically Noetherian representation. Then we show that there is an
ultrapower of this representation which has the same action type
quasi-identities as the original representation, but is not action type
geometrically equivalent to it.

\section{Some basics on modules and representations of groups.}

For the sake of completeness we will present in this section some
well-known basic notions and facts about the representations of
groups and modules which we will use later.

A representation $\left( V,G\right) $ is finitely generated if $G$
is a finitely generated group and $V$ is a finitely generated
$KG$-module.

The Cartesian product of the family of representations $\left\{ \left(
V_{i},G_{i}\right) \mid i\in I\right\} $ ($\left( V_{i},G_{i}\right) \in
Rep-K$ for every $i\in I$) is the representation $\left( \prod\limits_{i\in
I}V_{i},\prod\limits_{i\in I}G_{i}\right) $, with componentwise action.

If $\left\{ A_{i}\mid i\in I\right\} $ is a family of sets and $\mathfrak{F}$
is a filter over the set of indices $I$, then consider equivalence $\sim _{%
\mathfrak{F}}$ on the set $A=\prod\limits_{i\in I}A_{i}$: $a_{1}\sim _{%
\mathfrak{F}}a_{2}$ if $\left\{ i\in I\mid a_{1}^{\pi _{i}}=a_{2}^{\pi
_{i}}\right\} \in \mathfrak{F}$, where $a_{1},a_{2}\in A$, $\pi
_{i}:A\rightarrow A_{i}$ are projections. We denote by $\left[ a\right]
_{\sim _{\mathfrak{F}}}$ the  equivalence class $\sim _{\mathfrak{F}}$%
, generated by the element $a\in A$. The filtered product of the family of
sets $\left\{ A_{i}\mid i\in I\right\} $ by the filter $\mathfrak{F}$ is the
factor set $A/\sim _{\mathfrak{F}}=\left( \prod\limits_{i\in I}A_{i}\right)
/\sim _{\mathfrak{F}}$. If $\left\{ \left( V_{i},G_{i}\right) \mid i\in
I\right\} $ is a family of representations and $\mathfrak{F}$ is a filter
over the set of indices $I$, then the filtered product of this family of
representations by the filter $\mathfrak{F}$ is the representation%
\begin{equation*}
\left( \prod\limits_{i\in I}\left( V_{i},G_{i}\right) \right) /\sim _{%
\mathfrak{F}}=\left( \left( \prod\limits_{i\in I}V_{i}\right) /\sim _{%
\mathfrak{F}},\left( \prod\limits_{i\in I}G_{i}\right) /\sim _{\mathfrak{F}%
}\right) ,
\end{equation*}%
where action of the group $\left( \prod\limits_{i\in I}G_{i}\right) /\sim _{%
\mathfrak{F}}$ on the module $\left( \prod\limits_{i\in I}V_{i}\right) /\sim
_{\mathfrak{F}}$ is defined by $\left[ v\right] _{\sim _{\mathfrak{F}}}\circ %
\left[ g\right] _{\sim _{\mathfrak{F}}}=\left[ v\circ g\right] _{\sim _{%
\mathfrak{F}}}$ ($v\in \prod\limits_{i\in I}V_{i}$, $g\in \prod\limits_{i\in
I}G_{i}$).

The regular representation $\left( KG,G\right) $ is defined by: $v\circ g=vg$%
, where $v\in KG$, $g\in G$. If $U\leq KG_{KG}$ (in this way we denote both
a right ideal in a ring and a right submodule in a module) then $KG/U$ is
the $KG$-module and we have the representation $\left( KG/U,G\right) $: $%
v^{\nu }\circ g=\left( vg\right) ^{\nu }$, where $v\in $ $KG$, $g\in G$, $%
\nu :KG\rightarrow KG/U$ is the natural homomorphism.

Let $\rho :G\rightarrow \mathrm{Aut}_{K}\left( V\right) $ be the group
homomorphism, which defines a representation $\left( V,G\right) $. We denote%
\begin{equation*}
\ker \left( V,G\right) =\ker \rho =\left\{ g\in G\mid \forall v\in V\left(
v\circ g=v\right) \right\} .
\end{equation*}
The $\ker \left( V,G\right) $ is the normal subgroup of $G$. If $\ker \left(
V,G\right) =\left\{ 1\right\} $, then the representation is called faithful.
Let $\left( V,G\right) $ be an arbitrary representation, $\widetilde{G}%
=G/\ker \left( V,G\right) $, $\nu :G\rightarrow G/\ker \left( V,G\right) $
is the natural homomorphism. We can define the action of $\widetilde{G}$ on
module $V$ by $v\ast g^{\nu }=v\circ g$ ($v\in V$, $g\in G$). The
representation $\left( V,\widetilde{G}\right) $ is faithful and called the
faithful image of $\left(V,G\right) $.

Let $R$ be a ring with unit, $M$ is a $R$-module, $X\subseteq M$. We can
consider the annihilator of the set $X$: $\mathrm{ann}_{R}X=\left\{ r\in
R\mid \forall x\in X\left( xr=0\right) \right\} ,$ and the stabilizer of the
set $X$: $\mathrm{stab}_{R}X=\left\{ r\in R\mid \forall x\in X\left(
xr=x\right) \right\} $. It is well-known that $\mathrm{ann}_{R}X$ is the
right ideal of $R$ and if $X$ is a submodule of $M$, then the $\mathrm{ann}%
_{R}X$ is a two-sided ideal of $R$ and the $\mathrm{stab}_{R}X$ is a
semigroup of $R$. It is clear that
$\mathrm{stab}_{R}X=1+\mathrm{ann}_{R}X$.
If $\left( V,G\right) $ is a representation, $X\subseteq V$, we denote $%
\left( \mathrm{stab}_{KG}X\right) \cap G=\mathrm{stab}X$; by this notation
we have $\mathrm{stab}V=\ker \left( V,G\right) $.

\textbf{Proposition 1.1.} \textit{If }$R$\textit{\ is a ring with unit, }$%
U\leq R_{R}$\textit{\ right ideal of the ring }$R$\textit{, then }$\mathrm{%
ann}_{R}\left( R/U\right) $\textit{\ is the maximal two-sided ideal of the }$%
R$\textit{, contained in }$U$\textit{.} ([Pi], 2.1)

\textbf{Proposition 1.2.} \textit{Let }$K$\textit{\ be a commutative ring
with unit, }$G$\textit{\ be a group, }$U\leq KG_{KG}$\textit{, then }$\ker
\left( KG/U,G\right) $\textit{\ is the maximal normal subgroup of the }$G$%
\textit{\ contained in the group }$\left( 1+U\right) \cap G=\left\{ g\in
G\mid g-1\in U\right\} $\textit{.} This proposition is similar to \textbf{%
Proposition 1.1}.

\textbf{Corollary.} \textit{If }$U\leq _{KG}KG_{KG}$\textit{\ is two-sided
ideal of the }$KG$\textit{, then }%
\begin{equation*}
\left( 1+U\right) \cap G=\left\{ g\in G\mid g-1\in U\right\} =\ker \left(
KG/U,G\right) .
\end{equation*}%
If $\varphi :S\rightarrow R$ is a homomorphism of rings, then over every $R$%
-module $V_{R}$ we can define the structure of an $S$-module: $v\circ
s=vs^{\varphi }$ ($v\in V$, $s\in S$). We say in this case that $S$-module $%
V $ is defined by the homomorphism $\varphi $ and sometimes it is denoted by
$\left( V\right) _{\varphi }$.

\textbf{Proposition 1.3.} \textit{If }$\varphi :S\rightarrow R$\textit{\ is
an epimorphism of rings and }$U\leq R_{R}$\textit{\ is a right ideal of the
ring }$R$\textit{, then }$R/U\cong S/U^{\varphi ^{-1}}$\textit{\ as }$S$%
\textit{-modules. If }$U\leq _{R}R_{R}$\textit{\ a two-sided ideal, then }$%
R/U\cong S/U^{\varphi ^{-1}}$\textit{\ as rings.}

\textbf{Corollary 1.} \textit{Let }$\varphi :S\rightarrow R$\textit{\ is an
epimorphism of rings, }$U\leq R_{R}$\textit{\ a right ideal of the ring }$R$%
\textit{. Then}%
\begin{equation*}
S/\mathrm{ann}_{S}\left( R/U\right) \cong R/\left( \mathrm{ann}_{S}\left(
R/U\right) \right) ^{\varphi }=R/\mathrm{ann}_{R}\left( R/U\right)
\end{equation*}%
\textit{\ as rings.}

The epimorphism of groups $\varphi $ $:F\rightarrow G$ can be extended to
the epimorphism of associative algebras $\varphi :KF\rightarrow KG$. So, we
have

\textbf{Corollary 2.} \textit{If }$F=F\left( X\right) $\textit{\ is the free
}$n$\textit{-generated group, }$G$\textit{\ is another }$n$\textit{%
-generated group, }$U\leq KG_{KG}$\textit{\ is a right ideal of the ring }$%
KG $,\textit{\ then}%
\begin{equation*}
KF/\mathrm{ann}_{KF}\left( KG/U\right) \cong KG/\mathrm{ann}_{KG}\left(
KG/U\right)
\end{equation*}%
\textit{\ as associative algebras.}

\textbf{Proposition 1.4.} \textit{Let }$V_{R}\cong W_{R}$\textit{\ as }$R$%
\textit{-modules, then }$\mathrm{ann}_{R}\left( V_{R}\right) =\mathrm{ann}%
_{R}\left( W_{R}\right) $\textit{. }([Pi], 2.1)

\section{Operators on classes of representations of groups.}

Let $\mathfrak{X}$\ be a class of algebras in some variety $\Theta
$\ (many sorted in general). We consider the following operators on
classes of algebras:

$\mathcal{S}$: algebra $H\in \mathcal{S}\mathfrak{X}$ if and only if $H$\ is
a subalgebra of some algebra $G\in \mathfrak{X}$;

$\mathcal{C}$: algebra $H\in \mathcal{C}\mathfrak{X}$ if and only if $H$\ is
a Cartesian product of a family of algebras $\left\{ G_{i}\mid i\in
I\right\} $ from the class $\mathfrak{X}$;

$\mathcal{F}$: algebra $H\in \mathcal{F}\mathfrak{X}$ if and only if $H$\ is
a filtered product of a family of algebras $\left\{ G_{i}\mid i\in I\right\}
$ from the class $\mathfrak{X}$ by an arbitrary filter over the set $I$;

$\mathcal{C}_{up}$: algebra $H\in \mathcal{C}_{up}\mathfrak{X}$ if and only
if $H$\ is a filtered product of a family of algebras $\left\{ G_{i}\mid
i\in I\right\} $\ from the class $\mathfrak{X}$ by an arbitrary ultrafilter
over the set $I$;

$\mathcal{L}$: algebra $H\in \mathcal{L}\mathfrak{X}$ if and only if $%
H_{0}\in \mathfrak{X}$ for every finitely generated subalgebra $H_{0}\leq H$.

\textbf{Definitions 2.1. }\textit{Let }$\mathfrak{X}$\textit{\ be a class of
algebras in some variety }$\Theta $\textit{.}

\textit{If }$\mathcal{U}$\textit{\ is an operator on the classes of
algebras, we say that class of algebras }$\mathfrak{X}$\textit{\ is }\textbf{%
closed under the operator}\textit{\ }$\mathcal{U}$\textit{\ if }$\mathcal{U}%
\mathfrak{X}=\mathfrak{X}$\textit{.}

\textit{An operator }$\mathcal{U}$ \textit{on the classes of algebras is
called }\textbf{closed operator}\textit{\ if }$\mathcal{UU}\mathfrak{X}=%
\mathcal{U}\mathfrak{X}$\textit{\ for every class of algebras }$\mathfrak{X}$%
\textit{.}

\textit{An operator }$\mathcal{U}$ \textit{on the classes of algebras is
called }\textbf{monotone}\textit{\ if }$\mathcal{U}\mathfrak{X}_{1}\subset
\mathcal{U}\mathfrak{X}_{2}$\textit{\ holds when }$\mathfrak{X}_{1}\subset
\mathfrak{X}_{2}$\textit{\ (}$\mathfrak{X}_{1}$\textit{, }$\mathfrak{X}_{2}$%
\textit{\ - classes of algebras of the variety }$\Theta $\textit{).}

\textit{An operator }$\mathcal{U}$ \textit{on the classes of algebras is
called the }\textbf{operator of extension on }\textit{the fixed }\textbf{%
class} $\mathfrak{X}$\textit{\ if }$\mathcal{U}\mathfrak{X}\supset \mathfrak{%
X}$\textit{, an operator }$\mathcal{U}$ \textit{on the classes of algebras
is called the }\textbf{operator of extension}\textit{\ if it is an operator
of extension on all class of algebras.}

If $\mathcal{U}_{1},\ldots ,\mathcal{U}_{n}$ are operators on the classes of
algebras and $\mathfrak{X}$ is a class of algebras, we denote by $\left\{
\mathcal{U}_{1},\ldots ,\mathcal{U}_{n}\right\} \mathfrak{X}$ the minimal
class of algebras which contain the class $\mathfrak{X}$ and closed under
all operators $\mathcal{U}_{1},\ldots ,\mathcal{U}_{n}$. Of course, $\left\{
\mathcal{U}_{1},\ldots ,\mathcal{U}_{n}\right\} $ will also be the operator
on the classes of algebras.

It is clear that operators $\mathcal{L}$, $\mathcal{S}$, $\mathcal{C}$, $%
\mathcal{F}$, $\mathcal{C}_{up}$ are monotone. Operators $\mathcal{S}$, $%
\mathcal{C}$ are closed and\ operators of extension. $\mathcal{F}$ is also
an operator of extension, because over the set $\left\{ 1\right\} $ the
family of sets $\left\{ \left\{ 1\right\} \right\} $ is a filter. And $%
\mathcal{F}$ is a closed operator (see [Ma]).

It is well-known that for every class of algebras $\mathfrak{X}$ fulfills%
\begin{equation}
\mathcal{CS}\mathfrak{X}\subset \mathcal{SC}\mathfrak{X}  \tag{2.1}
\end{equation}%
and%
\begin{equation}
\mathcal{FS}\mathfrak{X}\subset \mathcal{SF}\mathfrak{X}.  \tag{2.2}
\end{equation}

About operator $\mathcal{L}$ in [PPT, Theorem 3] it was proved that if $%
\mathfrak{X}$ is a class of algebras, then%
\begin{equation}
\mathcal{SL}\mathfrak{X}=\mathcal{L}\mathfrak{X};  \tag{2.3}
\end{equation}%
\begin{equation}
\mathcal{CL}\mathfrak{X}\subset \mathcal{LSC}\mathfrak{X};  \tag{2.4}
\end{equation}%
\begin{equation}
\text{if }\mathcal{S}\mathfrak{X}=\mathfrak{X}\text{, then }\mathfrak{X}%
\subset \mathcal{L}\mathfrak{X}  \tag{2.5}
\end{equation}%
and%
\begin{equation}
\mathcal{LL}\mathfrak{X}=\mathcal{L}\mathfrak{X},  \tag{2.6}
\end{equation}%
and it was induced from this that $\left\{ \mathcal{L},\mathcal{S},\mathcal{C%
}\right\} =\mathcal{LSC}$.

Let $\mathfrak{X}$\ be a class of representations. On classes of
representations we can consider some special operators:

$\mathcal{Q}^{r}$: a representation $\left( V,G\right) \in \mathcal{Q}^{r}%
\mathfrak{X}$\ if and only if there exists a representation $\left(
V,D\right) \in \mathfrak{X}$, such that $\left( id_{V},\varphi \right)
:\left( V,D\right) \rightarrow \left( V,G\right) $ is a homomorphism of
representations\ and $\varphi :D\rightarrow G$\ is an epimorphism ([PV]
\textbf{1.3});

$\mathcal{Q}^{0}$: a representation $\left( V,G\right) \in \mathcal{Q}^{0}%
\mathfrak{X}$\ if and only if there exists a representation $\left(
V,D\right) \in \mathfrak{X}$, such that\textit{\ }$\left( id_{V},\varphi
\right) :\left( V,G\right) \rightarrow \left( V,D\right) $ is a homomorphism
of representations\ and $\varphi :G\rightarrow D$\ is an epimorphism ([PV]
\textbf{1.3});

$\mathcal{S}_{r}$: a representation $\left( V,G\right) \in \mathcal{S}_{r}%
\mathfrak{X}$\ if and only if $G\leq H$\ and $\left( V,H\right) \in
\mathfrak{X}$.

It is clear that operators $\mathcal{Q}^{0}$, $\mathcal{Q}^{r}$, $\mathcal{S}%
_{r}$ are monotone, closed and operators of extension.

Now we prove

\textbf{Lemma 2.1.}%
\begin{equation}
\mathcal{Q}^{r}\mathcal{L}\mathfrak{X}\subset \mathcal{LQ}^{r}\mathfrak{X},
\tag{2.7}
\end{equation}%
\begin{equation}
\mathcal{Q}^{0}\mathcal{L}\mathfrak{X}\subset \mathcal{LQ}^{0}\mathfrak{X},
\tag{2.8}
\end{equation}%
\begin{equation}
\mathcal{FQ}^{0}\mathfrak{X}\subset \mathcal{Q}^{0}\mathcal{F}\mathfrak{X},
\tag{2.9}
\end{equation}%
\begin{equation}
\mathcal{FQ}^{r}\mathfrak{X}\subset \mathcal{Q}^{r}\mathcal{F}\mathfrak{X},
\tag{2.10}
\end{equation}%
\textit{for every class of representations }$\mathfrak{X}$\textit{.}

\textit{Proof:}

Let $\mathfrak{X}$ be a class of representations.

Let $\left( V,G\right) \in \mathcal{Q}^{r}\mathcal{L}\mathfrak{X}$. Then
there exists representation $\left( V,D\right) \in \mathcal{L}\mathfrak{X}$
and epimorphism $\varphi :D\rightarrow G$, such that $\left( id_{V},\varphi
\right) $ is the homomorphism of representations. Let $\left(
V_{0},G_{0}\right) $ be a finitely generated subrepresentation of $\left(
V,G\right) $ and $G_{0}=\left\langle g_{1},\ldots ,g_{n}\right\rangle $. Let
$d_{i}^{\varphi }=g_{i}$ ($d_{i}\in D$, $1\leq i\leq n$). Denote $%
D_{0}=\left\langle d_{1},\ldots ,d_{n}\right\rangle $. $v\circ d=v\circ
d^{\varphi }$ for every $v\in V$, $d\in D$. So $\left( V_{0},D_{0}\right)
\leq \left( V,D\right) $ and $\left( V_{0},D_{0}\right) $ is a finitely
generated representation, because $V_{0}$ is finitely generated $KG_{0}$%
-module. $\left( V,D\right) \in \mathcal{L}\mathfrak{X}$, so $\left(
V_{0},D_{0}\right) \in \mathfrak{X}$. Therefore $\left( V_{0},G_{0}\right)
\in $ $\mathcal{Q}^{r}\mathfrak{X}$ and $\left( V,G\right) \in \mathcal{LQ}%
^{r}\mathfrak{X}$. (2.7) is proved.

Similarly we can prove (2.8).

Let $\left( V,G\right) \in \mathcal{FQ}^{0}\mathfrak{X}$. Then%
\begin{equation*}
\left( V,G\right) =\left( \left( \prod\limits_{i\in I}V_{i}\right) /\sim _{%
\mathfrak{F}},\left( \prod\limits_{i\in I}G_{i}\right) /\sim _{\mathfrak{F}%
}\right) ,
\end{equation*}%
where $\left\{ \left( V_{i},G_{i}\right) \mid i\in I\right\} \subset
\mathcal{Q}^{0}\mathfrak{X}$. So, an epimorphism $\varphi
_{i}:G_{i}\rightarrow D_{i}$ and a representation $\left( V_{i},D_{i}\right)
$, such that $\left( V_{i},D_{i}\right) \in \mathfrak{X}$ and $\left(
id_{V_{i}},\varphi _{i}\right) :\left( V_{i},G_{i}\right) \rightarrow \left(
V_{i},D_{i}\right) $ is a homomorphism of representations, exist for every $%
i\in I$. Hence,%
\begin{equation*}
\left( id_{\prod\limits_{i\in I}V_{i}},\varphi \right) :\left(
\prod\limits_{i\in I}V_{i},\prod\limits_{i\in I}G_{i}\right) \rightarrow
\left( \prod\limits_{i\in I}V_{i},\prod\limits_{i\in I}D_{i}\right) ,
\end{equation*}%
where $\varphi =\prod\limits_{i\in I}\varphi _{i}$, is a homomorphism of
representations. Let $\pi _{i}:\prod\limits_{i\in I}G_{i}\rightarrow G_{i}$
and $\rho _{i}:\prod\limits_{i\in I}D_{i}\rightarrow D_{i}$ be projections.
Let $g_{1},g_{2}\in \prod\limits_{i\in I}G_{i}$ and $g_{1}\sim _{\mathfrak{F}%
}g_{2}$, i.e. $\left\{ i\in I\mid g_{1}^{\pi _{i}}=g_{2}^{\pi _{i}}\right\}
\in \mathfrak{F}$. Then%
\begin{equation*}
\mathfrak{F\ni }\left\{ i\in I\mid g_{1}^{\varphi \rho _{i}}=g_{2}^{\varphi
\rho _{i}}\right\} \supset \left\{ i\in I\mid g_{1}^{\pi _{i}}=g_{2}^{\pi
_{i}}\right\} ,
\end{equation*}%
so $g_{1}^{\varphi }\sim _{\mathfrak{F}}g_{2}^{\varphi }$. Then we can define%
\begin{equation*}
\widetilde{\varphi }:\left( \prod\limits_{i\in I}G_{i}\right) /\sim _{%
\mathfrak{F}}\ni \left[ g\right] _{\sim _{\mathfrak{F}}}\rightarrow \left[
g^{\varphi }\right] _{\sim _{\mathfrak{F}}}\in \left( \prod\limits_{i\in
I}D_{i}\right) /\sim _{\mathfrak{F}}.
\end{equation*}%
$\widetilde{\varphi }$ is an epimorphism and $\left( id_{\left(
\prod\limits_{i\in I}V_{i}\right) /\sim _{\mathfrak{F}}},\widetilde{\varphi }%
\right) $ is a homomorphism of representations.%
\begin{equation*}
\left( \left( \prod\limits_{i\in I}V_{i}\right) /\sim _{\mathfrak{F}},\left(
\prod\limits_{i\in I}D_{i}\right) /\sim _{\mathfrak{F}}\right) \in \mathcal{F%
}\mathfrak{X},
\end{equation*}%
so%
\begin{equation*}
\left( V,G\right) =\left( \left( \prod\limits_{i\in I}V_{i}\right) /\sim _{%
\mathfrak{F}},\left( \prod\limits_{i\in I}G_{i}\right) /\sim _{\mathfrak{F}%
}\right) \in \mathcal{Q}^{0}\mathcal{F}\mathfrak{X}.
\end{equation*}%
(2.9) is proved.

Similarly we can prove (2.10).

\textbf{Corollary.}
\begin{equation}
\left\{ \mathcal{Q}^{r},\mathcal{Q}^{0},\mathcal{L},\mathcal{S},\mathcal{C}%
\right\} =\mathcal{LQ}^{0}\mathcal{Q}^{r}\mathcal{SC}  \tag{2.11}
\end{equation}%
\textit{and}%
\begin{equation}
\left\{ \mathcal{Q}^{0},\mathcal{Q}^{r},\mathcal{S},\mathcal{F}\right\} =%
\mathcal{Q}^{0}\mathcal{Q}^{r}\mathcal{SF}.  \tag{2.12}
\end{equation}

\textit{Proof:}

By the results of [PV] 1.3.2 and by (2.2), (2.9) and (2.10) we have
immediately (2.12).

To prove (2.11) we must use (2.3), (2.4), (2.5), (2.6), (2.7), (2.8), and
methods of [PPT, Theorem 3].

\section{Basic notions of algebraic geometry of representations.}

Let $\left( V,G\right) \in \Theta $ be a fixed representation. Let $\mathrm{%
Hom}\left( W_{\Theta }\left( X,Y\right) ,\left( V,G\right) \right) $ be the
affine space of finite rank over the $\left( V,G\right) $.

Consider in this affine space
the algebraic set%
\begin{equation*}
A=\left( T_{1},T_{2}\right) _{\left( V,G\right) }^{\prime }=\left\{ \left(
\alpha ,\beta \right) \in \mathrm{Hom}\left( W_{\Theta }\left( X,Y\right)
,\left( V,G\right) \right) \mid \ker \alpha \supset T_{1},\ker \beta \supset
T_{2}\right\} ,
\end{equation*}%
defined by an arbitrary pair $\left( T_{1},T_{2}\right) $ of sets of
equations ($T_{1}\subset E_{\Theta }(X,Y)$, $T_{2}\subset F_{\Theta }\left(
Y\right) $).

On the other hand, for an arbitrary set of points $A\subset \mathrm{Hom}%
\left( W_{\Theta }\left( X,Y\right) ,\left( V,G\right) \right) $ in the
affine space we have the "ideal of equations" (in our case - normal
subrepresentation of the $W_{\Theta }\left( X,Y\right) $ ), defined by this
set:
\begin{equation*}
T=A_{\left( V,G\right) }^{\prime }=\left( \bigcap\limits_{\left( \alpha
,\beta \right) \in A}\ker \alpha ,\bigcap\limits_{\left( \alpha ,\beta
\right) \in A}\ker \beta \right) .
\end{equation*}%
Then, the $\left( V,G\right) $-closure of the pair $\left(
T_{1},T_{2}\right) $ of sets of equations is:
\begin{equation*}
\left( T_{1},T_{2}\right) _{\left( V,G\right) }^{\prime \prime }=\left(
\bigcap\limits_{\left( \alpha ,\beta \right) \in \left( T_{1},T_{2}\right)
_{\left( V,G\right) }^{\prime }}\ker \alpha ,\bigcap\limits_{\left( \alpha
,\beta \right) \in \left( T_{1},T_{2}\right) _{\left( V,G\right) }^{\prime
}}\ker \beta \right) .
\end{equation*}%
The $\left( V,G\right) $-closure of the set of points $A\subset \mathrm{Hom}%
\left( W_{\Theta }\left( X,Y\right) ,\left( V,G\right) \right) $ is:%
\begin{equation*}
A_{(V,G)}^{\prime \prime }=\{(\mu ,\nu )\in \mathrm{Hom}(W_{\Theta
}(X,Y),(V,G))\mid \ker \mu \supset \bigcap\limits_{(\alpha ,\beta )\in
A}\ker \alpha ,
\end{equation*}%
\begin{equation*}
\ker \nu \supset \bigcap\limits_{(\alpha ,\beta )\in A}\ker \beta \}
\end{equation*}


We say that a pair of sets of equations $\left( S_{1},S_{2}\right) $ is
contained in a pair of sets of equations $\left( T_{1},T_{2}\right) $ ($%
S_{1},T_{1}\subset E_{\Theta }(X,Y)$, $S_{2},T_{2}\subset F_{\Theta }\left(
Y\right) $) and denote%
\begin{equation}
\left( S_{1},S_{2}\right) \subset \left( T_{1},T_{2}\right)  \tag{3.1}
\end{equation}%
if $S_{1}\subset T_{1}$ and $S_{2}\subset T_{2}$. The correspondence $%
^{\prime }$ is the Galois correspondence between sets of points and pairs of
sets of equations, that is:

1) $\left( \left( S_{1},S_{2}\right) \subset \left( T_{1},T_{2}\right)
\subset W_{\Theta }\left( X,Y\right) \right) \Rightarrow \left( \left(
S_{1},S_{2}\right) _{\left( V,G\right) }^{\prime }\supset \left(
T_{1},T_{2}\right) _{\left( V,G\right) }^{\prime }\right) $,

2) $\left( A\subset B\subset \mathrm{Hom}\left( W_{\Theta }\left( X,Y\right)
,\left( V,G\right) \right) \right) \Rightarrow \left( B_{\left( V,G\right)
}^{\prime }\supset A_{\left( V,G\right) }^{\prime }\right) $,

3) $\left( \left( T_{1},T_{2}\right) \subset W_{\Theta }\left( X,Y\right)
\right) \Rightarrow \left( \left( T_{1},T_{2}\right) \subset \left(
T_{1},T_{2}\right) _{\left( V,G\right) }^{\prime \prime }\right) $,

4) $\left( A\subset \mathrm{Hom}\left( W_{\Theta }\left( X,Y\right) ,\left(
V,G\right) \right) \right) \Rightarrow \left( A\subset A_{\left( V,G\right)
}^{\prime \prime }\right) $ -

for every $X$ and $Y$ and every representation $\left( V,G\right) $.

\textbf{Definitions 3.1.}\textit{\ We say that a pair of sets of equations }$%
\left( T_{1},T_{2}\right) $\textit{\ is }$\left( V,G\right) $\textbf{-closed}%
\textit{\ if }$\left( T_{1},T_{2}\right) =\left( T_{1},T_{2}\right) _{\left(
V,G\right) }^{\prime \prime }$\textit{\ and that a set of points } $A$%
\textit{\ is }$\left( V,G\right) $\textbf{-closed}\textit{\ if }$A=A_{\left(
V,G\right) }^{\prime \prime }$\textit{.}

If the pair of sets of equations $\left( T_{1},T_{2}\right) $\ is $\left(
V,G\right) $-closed, then $\left( T_{1},T_{2}\right) $ is the normal
subrepresentation of the $W_{\Theta }\left( X,Y\right) $.

As usual:

\textbf{Proposition 3.1. }The $\left( V,G\right) $-closure of a pair of sets
$\left( T_{1},T_{2}\right) $ is equal to the smallest $\left( V,G\right) $%
-closed pair containing the pair $\left( T_{1},T_{2}\right) $.

By [Pl2] (Proposition 3) we have

\textbf{Proposition 3.2.} \textit{Let }$\Theta _{1},\Theta _{2}$\textit{\
are a subvariety of }$Rep-K$\textit{, }$\left( V,G\right) \in \Theta
_{1}\subset \Theta _{2}$\textit{. There is a one-to-one order preserving
correspondence between lattices of }$\left( V,G\right) $\textit{-closed
subrepresentations in }$W_{\Theta _{2}}\left( X,Y\right) $\textit{\ and in }$%
W_{\Theta _{1}}(X,Y)$.

By this proposition we can consider the lattices of $\left( V,G\right) $%
-closed subrepresentations in the biggest variety of representations: in $%
Rep-K$.

Quasi-identity in $Rep-K$ can have the forms:%
\begin{equation}
\left( \left( \bigwedge\limits_{i=1}^{n_{1}}\left( w_{i}=0\right) \right)
\wedge \left( \bigwedge\limits_{i=1}^{n_{2}}\left( f_{i}=1\right) \right)
\right) \Rightarrow \left( w_{0}=0\right)  \tag{3.2.1}
\end{equation}%
and%
\begin{equation}
\left( \left( \bigwedge\limits_{i=1}^{n_{1}}\left( w_{i}=0\right) \right)
\wedge \left( \bigwedge\limits_{i=1}^{n_{2}}\left( f_{i}=1\right) \right)
\right) \Rightarrow \left( f_{0}=1\right) ,  \tag{3.2.2}
\end{equation}%
where $w_{i}\in XKF\left( Y\right) $ ($0\leq i\leq n_{1}$), $f_{i}\in
F\left( Y\right) $ ($0\leq i\leq n_{2}$). We say that a representation $%
\left( V,G\right) $ satisfies (3.2.1) and denote:%
\begin{equation*}
(V,G)\vDash \left( \left( \left( \bigwedge\limits_{i=1}^{n_{1}}\left(
w_{i}=0\right) \right) \wedge \left( \bigwedge\limits_{i=1}^{n_{2}}\left(
f_{i}=1\right) \right) \right) \Rightarrow \left( w_{0}=0\right) \right)
\end{equation*}%
%
%
%
%
%
%
%
%
%
%
%
%
%
%
\noindent if $w_{0}^{\alpha }=0$, for every homomorphism $\left( \alpha
,\beta \right) \in \mathrm{Hom}\left( W\left( X,Y\right) ,\left( V,G\right)
\right) $ which satisfies $w_{i}^{\alpha }=0$ for every $i\in \left\{
1,\ldots ,n_{1}\right\} $ and $f_{i}^{\beta }=1$ for every $i\in \left\{
1,\ldots ,n_{2}\right\} $. Similarly, a representation $\left( V,G\right) $
satisfies (3.2.2) if $f_{0}^{\beta }=1$, for every $\left( \alpha ,\beta
\right) \in \mathrm{Hom}\left( W\left( X,Y\right) ,\left( V,G\right) \right)
$ such that $w_{i}^{\alpha }=0$ for every $i\in \left\{ 1,\ldots
,n_{1}\right\} $ and $f_{i}^{\beta }=1$ for every $i\in \left\{ 1,\ldots
,n_{2}\right\} $.

Also we can consider the infinite "quasi-identities":%
\begin{equation}
\left( \left( \bigwedge\limits_{i\in I_{1}}\left( w_{i}=0\right) \right)
\wedge \left( \bigwedge\limits_{i\in I_{2}}\left( f_{i}=1\right) \right)
\right) \Rightarrow \left( w_{0}=0\right)  \tag{3.2.1'}
\end{equation}%
and%
\begin{equation}
\left( \left( \bigwedge\limits_{i\in I_{1}}\left( w_{i}=0\right) \right)
\wedge \left( \bigwedge\limits_{i\in I_{2}}\left( f_{i}=1\right) \right)
\right) \Rightarrow \left( f_{0}=1\right) ,  \tag{3.2.2'}
\end{equation}%
where $w_{i}\in XKF\left( Y\right) $ ($i\in I_{1}\cup \left\{ 0\right\} $), $%
f_{i}\in F\left( Y\right) $ ($i\in I_{2}\cup \left\{ 0\right\} $), $%
I_{1},I_{2}$ is not necessary finite. We say that a representation $\left(
V,G\right) $ satisfies (3.2.1') if $w_{0}^{\alpha }=0$, for every $\left(
\alpha ,\beta \right) \in \mathrm{Hom}\left( W\left( X,Y\right) ,\left(
V,G\right) \right) $ such that $w_{i}^{\alpha }=0$ for every $i\in I_{1}$
and $f_{i}^{\beta }=1$ for every $i\in I_{2}$. Similarly, a representation $%
\left( V,G\right) $ satisfies (3.2.2') if $f_{0}^{\beta }=1$, for every $%
\left( \alpha ,\beta \right) \in \mathrm{Hom}\left( W\left( X,Y\right)
,\left( V,G\right) \right) $, such that $w_{i}^{\alpha }=0$ for every $i\in
I_{1}$ and $f_{i}^{\beta }=1$ for every $i\in I_{2}$.

Let $\mathfrak{X}$ be a class of representations. Denote by $qId\mathfrak{X}$
the set of quasi-identities satisfied by all representations of this class.
Let now $\mathfrak{Q}$ be a set of quasi-identities of representations.
Denote by $qVar\mathfrak{Q}$ the class of all representations, which satisfy
all quasi-identities from the $\mathfrak{Q}$; this class is called a
quasi-variety of representations. The quasi-variety generated by the class $%
\mathfrak{X}$ 
is denoted by $qVar\mathfrak{X}$. If $\mathfrak{X}\mathcal{=}\left\{ \left(
V,G\right) \right\} $ then we denote $qId\mathfrak{X}=qId\left( V,G\right) $%
, $qVar\mathfrak{X}=qVar\left( V,G\right) $.

It is easy to see that

\textbf{Proposition 3.3.}\textit{\ If }$\left\{ w_{i}\mid i\in I_{1}\right\}
\cup \left\{ w_{0}\right\} \subset XKF\left( Y\right) $\textit{, }$\left\{
f_{i}\mid i\in I_{2}\right\} \cup \left\{ f_{0}\right\} \subset F\left(
Y\right) $\textit{, }$\left( T_{1},T_{2}\right) =\left( \left\{ w_{i}\mid
i\in I_{1}\right\} ,\left\{ f_{i}\mid i\in I_{2}\right\} \right)
_{(V,G)}^{\prime \prime }$\textit{\ then }
\begin{equation*}
(V,G)\vDash \left( \left( \left( \bigwedge\limits_{i\in I_{1}}\left(
w_{i}=0\right) \right) \wedge \left( \bigwedge\limits_{i\in I_{2}}\left(
f_{i}=1\right) \right) \right) \Rightarrow \left( w_{0}=0\right) \right)
\end{equation*}%
%
%
%
%
%
%
%
%
%
%
%
%
%
%
\noindent \textit{if and only if }$w_{0}\in T_{1}$\textit{\ and}%
\begin{equation*}
(V,G)\vDash \left( \left( \left( \bigwedge\limits_{i\in I_{1}}\left(
w_{i}=0\right) \right) \wedge \left( \bigwedge\limits_{i\in I_{2}}\left(
f_{i}=1\right) \right) \right) \Rightarrow \left( f_{0}=1\right) \right)
\end{equation*}%
%
%
%
%
%
%
%
%
%
%
%
%
%
%
\noindent \textit{if and only if }$f_{0}\in T_{2}$\textit{.}%


\textbf{Proposition 3.4.} $\left( T_{1},T_{2}\right) \subset W\left(
X,Y\right) $\textit{\ is a }$\left( V,G\right)$-closed representation if and
only if $\left( XKF\left( Y\right) /T_{1},F\left( Y\right) /T_{2}\right)
=W\left( X,Y\right) /\left( T_{1},T_{2}\right) \in \mathcal{SC}\left(
V,G\right) $\textit{.}

\textit{Proof:}

We apply the Remake theorem for representations to the representation $%
W\left( X,Y\right) /\left( T_{1},T_{2}\right) $.

\textbf{Corollary.} \textit{Let }$\left( V,G\right) ,\left( W,H\right) \in
Rep-K$\textit{. Every }$\left( W,H\right) $\textit{-closed representation is
a }$\left( V,G\right) $\textit{-closed representation if and only if }$%
\left( W,H\right) \in \mathcal{LSC}\left( V,G\right) $\textit{.}

\textit{Proof:}

Let $\left( W,H\right) \in \mathcal{LSC}\left( V,G\right) $. Let $\left(
T_{1},T_{2}\right) \subset W\left( X,Y\right) $\ be a $\left( W,H\right) $%
-closed representation, then $W\left( X,Y\right) /\left( T_{1},T_{2}\right)
\in \mathcal{SC}\left( W,H\right) $. $\left\{ \mathcal{L},\mathcal{S},%
\mathcal{C}\right\} =\mathcal{LSC}$ ([PPT], Theorem 3), so $W\left(
X,Y\right) /\left( T_{1},T_{2}\right) \in \mathcal{LSC}\left( V,G\right) $.
The representation $W\left( X,Y\right) /\left( T_{1},T_{2}\right) $ is
finitely generated, so it belongs to $\mathcal{SC}\left( V,G\right) $ and $%
\left( T_{1},T_{2}\right) $ is a $\left( V,G\right) $-closed representation.

Let every $\left( W,H\right) $-closed representation be a $\left( V,G\right)
$-closed representation. Let $\left( W_{0},H_{0}\right) \subset \left(
W,H\right) $ be a finitely generated subrepresentation. There are $X$ and $Y$%
, such that $\left( W_{0},H_{0}\right) \cong W\left( X,Y\right) /\left(
T_{1},T_{2}\right) $, where $\left( T_{1},T_{2}\right) $ is a normal
subrepresentation of $W\left( X,Y\right) $. By \textbf{Proposition 3.4 }$%
\left( T_{1},T_{2}\right) $ is a $\left( W,H\right) $-closed. So $\left(
T_{1},T_{2}\right) $ is a $\left( V,G\right) $-closed. Therefore, $\left(
W_{0},H_{0}\right) \cong \left( XKF\left( Y\right) /T_{1},F\left( Y\right)
/T_{2}\right) \in \mathcal{SC}\left( V,G\right) $ and $\left( W,H\right) \in
\mathcal{LSC}\left( V,G\right) $. The proof is complete.

\textbf{Definition 3.2.} \textit{Representations }$\left( V_{1},G_{1}\right)
$\textit{\ and }$\left( V_{2},G_{2}\right) $\textit{\ are called }\textbf{%
geometrically equivalent}\textit{\ (}denoted\textit{\ }$\left(
V_{1},G_{1}\right) \sim \left( V_{2},G_{2}\right) $\textit{) if }$\left(
T_{1},T_{2}\right) _{\left( V_{1},G_{1}\right) }^{\prime \prime }=\left(
T_{1},T_{2}\right) _{\left( V_{2},G_{2}\right) }^{\prime \prime }$\textit{\
for every }$X$\textit{\ and }$Y$\textit{\ and for every pair of sets }$%
\left( T_{1},T_{2}\right) \subset W\left( X,Y\right) $\textit{.}

\textbf{Corollary} from \textbf{Proposition 3.1.} $\left( V_{1},G_{1}\right)
\sim \left( V_{2},G_{2}\right) $\textit{\ if and only if every }$\left(
V_{1},G_{1}\right) $\textit{-closed representation is a }$\left(
V_{2},G_{2}\right) $\textit{-closed representation and vice versa.}

This corollary and \textbf{Proposition 3.2} imply that if two
representations $\left( V_{1},G_{1}\right) $ and $\left( V_{2},G_{2}\right) $
belong to the same variety $\Theta $, then their geometrical equivalence
does not depend on the fact in which affine space: absolute
or relative 
we consider algebraic geometries over these representations.

\textbf{Corollary 1} from \textbf{Proposition 3.3.} \textit{If }$(
V_{1},G_{1}) \sim ( V_{2},G_{2}) $\textit{\ then }$qId( V_{1},G_{1})$ $=qId(
V_{2},G_{2}) $\textit{.}

\textbf{Corollary 2} from \textbf{Proposition 3.3.} $\left(
V_{1},G_{1}\right) \sim \left( V_{2},G_{2}\right) $\textit{\ if and only if
representations }$\left( V_{1},G_{1}\right) $\textit{\ and }$\left(
V_{2},G_{2}\right) $\textit{\ have same infinite "quasi-identities".}

By [PPT, Theorem 3], we have:

\textbf{Proposition 3.5.} \textit{Let }$\left( V_{1},G_{1}\right) $\textit{\
and }$\left( V_{2},G_{2}\right) $\textit{\ be representations. Then }$\left(
V_{1},G_{1}\right) \sim \left( V_{2},G_{2}\right) $\textit{\ if and only if }%
$\mathcal{LSC}\left( V_{1},G_{1}\right) =\mathcal{LSC}\left(
V_{2},G_{2}\right) $\textit{.}

\textbf{Corollary.} \textit{Let }$\left( V_{1},G_{1}\right) $\textit{\ and }$%
\left( V_{2},G_{2}\right) $\textit{\ be representations. If }$\left(
V_{1},G_{1}\right) \sim \left( V_{2},G_{2}\right) $\textit{\ then }$%
G_{1}\sim G_{2}$\textit{\ as groups.}

\textit{Proof:} Let $G_{1}^{0}\leq G_{1}$ be a finitely generated subgroup.
There is a finitely generated subrepresentation $\left(
V_{1}^{0},G_{1}^{0}\right) \subset \left( V_{1},G_{1}\right) $ and there is
an embedding of representations $\left( V_{1}^{0},G_{1}^{0}\right)
\hookrightarrow \left( V_{2},G_{2}\right) ^{I}=\left(
V_{2}^{I},G_{2}^{I}\right) $ ($I$ - some set of indices). So, there is the
embedding of groups: $G_{1}^{0}\hookrightarrow G_{2}^{I}$. Therefore $%
G_{1}\in \mathcal{LSC}\left( G_{2}\right) $. By symmetry, $G_{2}\in \mathcal{%
LSC}\left( G_{1}\right) $. By [PPT, Theorem 3], the proof is complete.

Define now the notions of a Noetherian variety of representations
(subvariety of $Rep-K$), geometrically Noetherian representation and
logically Noetherian representation:

\textbf{Definitions 3.3. }\textit{We call a variety }$\Theta \subset Rep-K$%
\textit{\ }\textbf{Noetherian}\textit{\ if for every }$X$\textit{\ \ and }$Y$%
\textit{\ every normal subrepresentation of }$W_{\Theta }(X,Y)$\textit{\ is
finitely generated.}

\textit{A representation }$\left( V,G\right) \in \Theta $ \textit{is called }%
\textbf{geometrically Noetherian}\textit{\ if for every sets }$X$\textit{\
and }$Y$\textit{\ and every pair of sets }$\left( T_{1},T_{2}\right) \subset
W_{\Theta }\left( X,Y\right) $\textit{, there is a pair of finite subsets }$%
\left( R_{1},R_{2}\right) \subset \left( T_{1},T_{2}\right) $\textit{, such
that }$\left( T_{1},T_{2}\right) _{\left( V,G\right) }^{\prime }=\left(
R_{1},R_{2}\right) _{\left( V,G\right) }^{\prime }$\textit{.}

\textit{A representation }$\left( V,G\right) \in \Theta $\textit{\ is called
}\textbf{logically Noetherian}\textit{\ if for every sets }$X$\textit{\ and }%
$Y$\textit{,\ every pair of sets }$\left( T_{1},T_{2}\right) \subset
W_{\Theta }\left( X,Y\right) $\textit{\ and every }$w\in E_{\Theta }\left(
X,Y\right) $\textit{\ (}$f\in F_{\Theta }\left( Y\right) $\textit{) belongs
to the first (second) component of the pair }$\left( T_{1},T_{2}\right)
_{\left( V,G\right) }^{\prime \prime }$\textit{\ exists a pair of finite
subsets }$\left( R_{1},R_{2}\right) \subset \left( T_{1},T_{2}\right) $%
\textit{, such that }$w$\textit{\ (}$f$\textit{) belongs to the first
(second) component of the pair }$\left( R_{1},R_{2}\right) _{\left(
V,G\right) }^{\prime \prime }$\textit{.}

It is clear that Noetherianity of the variety $\Theta $ of representations
is equivalent to the ascending chain condition for normal subrepresentations
in every finitely generated relatively free representation $W_{\Theta }(X,Y)$
and geometrical Noetherianity of the representation $\left( V,G\right) $ is
equivalent to the ascending chain condition for $\left( V,G\right) $-closed
normal subrepresentations in every finitely generated free representation $%
W_{\Theta }\left( X,Y\right) $ .

By (3.1) the order on a family of pairs of sets $\left\{ \left( R,T\right)
\mid \left( R,T\right) \subset W\left( X,Y\right) \right\} $ is defined. We
can consider directed systems by this order. Also we can consider the union
of two pairs of sets: $\left( R_{1},T_{1}\right) \cup \left(
R_{2},T_{2}\right) =\left( R_{1}\cup R_{2},T_{1}\cup T_{2}\right) $.
According to [Pl4, Proposition 7]:

\textbf{Proposition 3.6.} \textit{A representation }$\left( V,G\right) \in
\Theta $\textit{\ is logically Noetherian if and only if the union of any
directed system of }$\left( V,G\right) $\textit{-closed subrepresentations
in the }$W_{\Theta }\left( X,Y\right) $\textit{\ for every }$X$ and $Y$%
\textit{\ is also a }$\left( V,G\right) $\textit{-closed subrepresentation.}

So, by \textbf{Proposition 3.2}, geometric Noetherianity and logic
Noetherianity of representation is not depend in what subvariety $\Theta $
we consider those. Also every representation $\left( V,G\right) $ from the
Noetherian variety $\Theta $ is geometrically Noetherian. And every
geometrically Noetherian representation is logically Noetherian.

By \textbf{Proposition 3.3}, if a representation $\left( V,G\right) $ is
logically Noetherian then for every infinite quasi-identity of the form
(3.2.1') ( (3.2.2') ) which is fulfilled in the $\left( V,G\right) $ there
is a finite quasi-identity with the minor premise and the same conclusion
which is also fulfilled in the $\left( V,G\right) $. A representation $%
\left( V,G\right) $ is geometrically Noetherian if and only if in this
reduction choosing of a premise is not dependent on the conclusion, but only
on the infinite premise.

\textbf{Proposition 3.7.} \textit{If a representation }$\left( V,G\right) $
\textit{is geometrically (logically) Noetherian, then the group }$G$\textit{%
\ is geometrically (logically) Noetherian too.}

\textit{Proof:}

Let $\left( V,G\right) $ be a logically Noetherian representation. Let the
group $G$ satisfies the infinite group quasi-identity%
\begin{equation*}
\left( \bigwedge\limits_{i\in I}\left( f_{i}=1\right) \right) \Rightarrow
\left( f_{0}=1\right) ,
\end{equation*}%
where $\left\{ f_{i}\mid i\in I\right\} \cup \left\{ f_{0}\right\} \subset
F\left( Y\right) $. This quasi-identity can be considered as a special case
of (3.2.2') and can be reduced in $\left( V,G\right) $ to the finite
quasi-identity%
\begin{equation*}
\left( \bigwedge\limits_{i\in I_{0}}\left( f_{i}=1\right) \right)
\Rightarrow \left( f_{0}=1\right) ,
\end{equation*}%
where $I_{0}\subset I$, $\left\vert I_{0}\right\vert <\aleph _{0}$. Every
group homomorphism $\beta \in \mathrm{Hom}\left( F\left( Y\right) ,G\right) $
can be realized as the second component in a homomorphism of representations
$\left( \alpha ,\beta \right) \in \mathrm{Hom}\left( W\left( X,Y\right)
,\left( V,G\right) \right) $, for example, as the second component in the
pair $\left( 0,\beta \right) \in \mathrm{Hom}\left( W\left( X,Y\right)
,\left( V,G\right) \right) $ (for every $X$). Therefore, $G$ fulfill the
quasi-identity%
\begin{equation*}
\left( \bigwedge\limits_{i\in I_{0}}\left( f_{i}=1\right) \right)
\Rightarrow \left( f_{0}=1\right) .
\end{equation*}%
Hence, the group $G$ is logically Noetherian.

Analogously, we can prove that if representation $\left( V,G\right) $ is
geometrically Noetherian, then the group $G$\ is geometrically Noetherian
too. The proof is complete.

This proposition and the \textbf{Corollary} from \textbf{Proposition 3.5}
show that algebraic geometry over representations of groups in the regular
sense, i.e. the algebraic geometry which deals with equations on acting
groups and action-type equations, is very closely connected with the
algebraic geometry over groups. For example, if a group $G$ is non
geometrically (logically) Noetherian, then every representation of this
group is non geometrically (logically) Noetherian and this fact does not
depend on the action of this group on a module. So, in order to study the
geometry which enjoys the peculiarities of the action one has to consider
not the "two-sided" geometry above, but the one-sided action-type geometry.


\section{Basic notions of action type algebraic geometry of representations.}

In action type algebraic geometry of representations, we consider algebraic
sets in the affine space $\mathrm{Hom}\left( W_{\Theta }\left( X,Y\right)
,\left( V,G\right) \right) $ defined only by action type equations: $w=0$ -
where $w\in E_{\Theta }(X,Y)$.

We have, as above, the Galois correspondence between sets of "points" and
sets of action type equations:

1) $\left( S\subset T\subset E_{\Theta }(X,Y)\right) \Rightarrow \left(
S_{\left( V,G\right) }^{\nabla }\supset T_{\left( V,G\right) }^{\nabla
}\right) $,

2) $\left( A\subset B\subset \mathrm{Hom}\left( W_{\Theta }\left( X,Y\right)
,\left( V,G\right) \right) \right) \Rightarrow \left( B_{\left( V,G\right)
}^{\nabla }\supset A_{\left( V,G\right) }^{\nabla }\right) $,

3) $\left( T\subset E_{\Theta }(X,Y)\right) \Rightarrow \left( T\subset
T_{\left( V,G\right) }^{\nabla \nabla }\right) $,

4) $\left( A\subset \mathrm{Hom}\left( W_{\Theta }\left( X,Y\right) ,\left(
V,G\right) \right) \right) \Rightarrow \left( A\subset A_{\left( V,G\right)
}^{\nabla \nabla }\right) $ -

for every $X$ and $Y$ and every arbitrary representation $\left( V,G\right)
\in \Theta $.

\textbf{Definition 4.1.}\textit{\ We say that a set of action type equations
}$T$\textit{\ is }\textbf{action type }$\left( V,G\right) $\textbf{-closed}%
\textit{\ if }$T=T_{\left( V,G\right) }^{\nabla \nabla }$\textit{.}

If the set of equations $T$\ is action type $\left( V,G\right) $-closed,
then $T$ is a $KF_{\Theta }\left( Y\right) $-submodule of the $E_{\Theta
}(X,Y)$.

The Galois correspondence implies

\textbf{Proposition 4.1. }\textit{The action type }$\left( V,G\right) $%
\textit{-closure of the set }$T$\textit{\ is equal to the smallest action
type }$\left( V,G\right) $\textit{-closed submodule, containing the set }$T$%
\textit{.}

\textbf{Proposition 4.2.} \textit{Let }$\left( V,G\right) \in \Theta $%
\textit{. A }$KF_{\Theta }\left( Y\right) $\textit{-submodule }$T\leq
E_{\Theta }(X,Y)$\textit{\ is an action type }$\left( V,G\right) $\textit{%
-closed\ if and only if there exists a normal subgroup }$H\trianglelefteq
F_{\Theta }\left( Y\right) $\textit{\ such that }$\left( T,H\right) \subset
W_{\Theta }\left( X,Y\right) $\textit{\ is the }$\left( V,G\right) $\textit{%
-closed subrepresentation.}

\textit{Proof:}

It is clear that $\left( T,H\right) _{\left( V,G\right) }^{\prime }\subset
T_{\left( V,G\right) }^{\nabla }$ for every $H\subset F_{\Theta }\left(
Y\right) $. So, if $\left( T,H\right) =\left( T,H\right) _{\left( V,G\right)
}^{\prime \prime }$ then $T\supset T_{\left( V,G\right) }^{\nabla \nabla
}\supset T$. Therefore $T$ is an action type $\left( V,G\right) $-closed
submodule.

Let $T=T_{\left( V,G\right) }^{\nabla \nabla }$. Denote $\bigcap\limits_{%
\beta \in \mathrm{Hom}\left( F_{\Theta }\left( Y\right) ,G\right) }\ker
\beta =Id_{\Theta }\left( G,Y\right) $. It is clear that $\left(
T,Id_{\Theta }\left( G,Y\right) \right) _{\left( V,G\right) }^{\prime
}=T_{\left( V,G\right) }^{\nabla }$. Thus, $\left( T,Id_{\Theta }\left(
G,Y\right) \right) _{\left( V,G\right) }^{\prime \prime }=\left(
T,Id_{\Theta }\left( G,Y\right) \right) $ since $\left( 0,\beta \right) \in
T_{\left( V,G\right) }^{\nabla }$ for every group homomorphism $\beta
:F_{\Theta }\left( Y\right) \rightarrow G$. The proof is complete.

\textbf{Corollary 1.} $T\leq XKF\left( Y\right) $\textit{\ is an action type
}$\left( V,G\right) $\textit{-closed submodule if and only if there exists a
normal subgroup }$H\trianglelefteq F\left( Y\right) $\textit{, such that }$%
\left( XKF\left( Y\right) /T,F\left( Y\right) /H\right) \in \mathcal{SC}%
\left( V,G\right) $\textit{.}

\textit{Proof:} By \textbf{Proposition 3.4}.

\textbf{Remark 4.1.} We can see from the proof, that in \textbf{Proposition
4.2} and its \textbf{Corollary 1} one can always take $\bigcap\limits_{\beta
\in \mathrm{Hom}\left( F_{\Theta }\left( Y\right) ,G\right) }\ker \beta
=Id_{\Theta }\left( G,Y\right) $ as normal subgroup $H\trianglelefteq
F_{\Theta }\left( Y\right) $ .

From \textbf{Proposition 4.2} and \textbf{Proposition 3.2} we can easy
conclude

\textbf{Proposition 4.3.} \textit{Let }$\Theta _{1},\Theta _{2}$\textit{\
are a subvariety of }$Rep-K$\textit{, }$\left( V,G\right) \in \Theta
_{1}\subset \Theta _{2}$\textit{. There is a one-to-one order preserving
correspondence between lattices of action type }$\left( V,G\right) $\textit{%
-closed submodules in }$E_{\Theta _{2}}\left( X,Y\right) $\textit{\ and in }$%
E_{\Theta _{1}}(X,Y)$.

By this proposition we can consider the lattices of action type $\left(
V,G\right) $-closed submodules in the baggiest variety of representations:
in $Rep-K$.

We have immediately

\textbf{Proposition 4.4.}
\begin{equation*}
\left( V,G\right) \vDash \left( \left( \bigwedge\limits_{i\in I}\left(
w_{i}=0\right) \right) \Rightarrow \left( w_{0}=0\right) \right) ,
\end{equation*}%
\noindent \textit{\ where }$\left\{ w_{0}\right\} \cup \left\{ w_{i}\mid
i\in I\right\} \subset XKF\left( Y\right) $\textit{, if and only if }$%
w_{0}\in \left\{ w_{i}\mid i\in I\right\} _{\left( V,G\right) }^{\nabla
\nabla }$\textit{.}

\textbf{Definition 4.2.} \textit{We say that a quasi-variety of
representations }$\mathfrak{X}$\textit{\ is an }\textbf{action type
quasi-variety}\textit{\ if it can be defined by a set of action type
quasi-identities.}

It means that $\mathfrak{X}$\textit{\ }is an action type quasi-variety of
representations\textit{\ }if and only if there exists a set of action type
quasi-identities $\mathfrak{Q}$ such that $\mathfrak{X}=qVar\mathfrak{Q}$.

The set of all action type quasi-identities satisfied by representation $%
\left( V,G\right) $ is denoted by $qId_{a.t.}\left( V,G\right) $. Let $%
\mathfrak{X}\subset Rep-K$ be a class of representations. Denote by $%
qId_{a.t.}\mathfrak{X}$ the set of action type quasi-identities satisfied by
all representations from $\mathfrak{X}$. Clearly, $qId_{a.t.}\mathfrak{X}%
=\dbigcap\limits_{\left( V,G\right) \in \mathfrak{X}}qId_{a.t.}\left(
V,G\right) $ . We denote $qVar\left( qId_{a.t.}\mathfrak{X}\right)
=qVar_{a.t.}\mathfrak{X}$.

\textbf{Definition 4.3.} \textit{The action type quasi-variety }$qVar_{a.t.}%
\mathfrak{X}$\textit{\ we call }\textbf{action type quasi-variety, generated
by the class}\textit{\ }$\mathfrak{X}$\textit{.}

\section{Action type geometrical equivalence of representations.}

\textbf{Corollary} from \textbf{Proposition 4.1.} \textit{Let }$\left(
V_{1},G_{1}\right) ,\left( V_{2},G_{2}\right) \in \Theta $\textit{. Then }$%
\left( V_{1},G_{1}\right) \sim _{a.t.}\left( V_{2},G_{2}\right) $\textit{\
if and only if for every finite }$X$\textit{\ and }$Y$\textit{\ every action
type }$\left( V_{1},G_{1}\right) $\textit{-closed submodule of }$E\left(
X,Y\right) $\textit{\ is the action type }$\left( V_{2},G_{2}\right) $%
\textit{-closed submodule and vice versa. }

By this \textbf{Corollary} and by \textbf{Proposition 4.3}, action type
geometrical equivalence of representations $\left( V_{1},G_{1}\right) $ and $%
\left( V_{2},G_{2}\right) $ can be recognized in all subvariety $\Theta
\subseteq Rep-K$, such that $\left( V_{1},G_{1}\right) ,(V_{2},G_{2})\in
\Theta $. Below we use for this purpose the biggest variety of
representations: $Rep-K$.

\textbf{Corollary 1} from \textbf{Proposition 4.4.} \textit{If }$\left(
V_{1},G_{1}\right) \sim _{a.t.}\left( V_{2},G_{2}\right) $\textit{\ then }$%
qId_{a.t.}\left( V_{1},G_{1}\right) =qId_{a.t.}\left( V_{2},G_{2}\right) $%
\textit{.}

\textbf{Corollary 2} from \textbf{Proposition 4.4.} $\left(
V_{1},G_{1}\right) \sim _{a.t.}\left( V_{2},G_{2}\right) $\textit{\ if and
only if representations }$\left( V_{1},G_{1}\right) $\textit{\ and }$\left(
V_{2},G_{2}\right) $\textit{\ have the same infinite action type
quasi-identities.}

Also we have

\textbf{Corollary 2 }from \textbf{Proposition 4.2.} \textit{If two
representations }$\left( V_{1},G_{1}\right) $\textit{\ and }$\left(
V_{2},G_{2}\right) $\textit{\ are geometrically equivalent then they are
action type\ geometrically equivalent.}

\textbf{Remark 5.1.} In spite of this \textbf{Corollary} and \textbf{%
Corollary} from \textbf{Proposition 3.5}, if two representations $\left(
V_{1},G_{1}\right) $ and $\left( V_{2},G_{2}\right) $ are action type\
geometrically equivalent and groups $G_{1}$ and $G_{2}$ are geometrically
equivalent, the representations $\left( V_{1},G_{1}\right) $ and $\left(
V_{2},G_{2}\right) $ are not necessarily geometrically equivalent.

\textbf{Definitions 5.1.}\textit{\ Two representations }$\left(
V_{1},G_{1}\right) $\textit{\ and }$\left( V_{2},G_{2}\right) $\textit{\ are
called }\textbf{(algebraically) equivalent}\textit{\ if the corresponding
faithful representations are isomorphic.}

\textit{A class of representations }$\mathfrak{X}$\textit{\ is called }%
\textbf{saturated}\textit{\ if with a representation }$\left( V,G\right) \in
\mathfrak{X}$\textit{\ it contains all representations which are
algebraically equivalent to the representation }$\left( V,G\right) $\textit{.%
}

By [PV, \textbf{1.3}], two representations $\left( V_{1},G_{1}\right) $\ and
$\left( V_{2},G_{2}\right) $\ are algebraically equivalent\ if and only if $%
\mathcal{Q}^{0}\mathcal{Q}^{r}\left( V_{1},G_{1}\right) =\mathcal{Q}^{0}%
\mathcal{Q}^{r}\left( V_{2},G_{2}\right) $, so a class of representations $%
\mathfrak{X}$ is saturated if and only if $\left\{ \mathcal{Q}^{0},\mathcal{Q%
}^{r}\right\} \mathfrak{X=\mathcal{Q}}^{0}\mathcal{Q}^{r}\mathfrak{X=X}$.

\textbf{Theorem 5.1.} \textit{Let }$\left( Z,H\right) ,\left( V,G\right) \in
Rep-K$\textit{. Every action type }$\left( Z,H\right) $\textit{-closed
submodule }$T\leq \left( XKF\left( Y\right) \right) _{KF\left( Y\right) }$
is an action type $\left( V,G\right) $\textit{-closed submodule if and only
if }$\left( Z,H\right) \in \mathcal{LQ}^{0}\mathcal{Q}^{r}\mathcal{SC}\left(
V,G\right) $.

\textit{Proof:}

Let $\left( Z,H\right) \in \mathcal{LQ}^{0}\mathcal{Q}^{r}\mathcal{SC}\left(
V,G\right) $. Let $T\leq \left( XKF\left( Y\right) \right) _{KF\left(
Y\right) }$ be an action type $\left( Z,H\right) $-closed submodule. By
\textbf{Proposition 4.2}, there exists a normal subgroup $P\trianglelefteq
F\left( Y\right) $ such that $\left( T,P\right) \subset W\left( X,Y\right) $%
\ is the $\left( Z,H\right) $-closed subrepresentation. By \textbf{%
Proposition 3.4} and by (2.11) $(XKF(Y)/T,F(Y)/P)$$\in \mathcal{SC}%
(Z,H)\subset \mathcal{LQ}^{0}\mathcal{Q}^{r}\mathcal{SC}(V,G)$. $%
(XKF(Y)/T,F(Y)/P)$ is finitely generated, so, by projectivity of the free
groups, there exists $S\trianglelefteq F\left( Y\right) $ such that $\left(
XKF\left( Y\right) /T,F\left( Y\right) /S\right) \in \mathcal{SC}\left(
V,G\right) $. So, by \textbf{Proposition 3.4} and \textbf{Proposition 4.2}, $%
T$ is an action type $\left( V,G\right) $-closed submodule.

Let every action type $\left( Z,H\right) $-closed submodule be an action
type $\left( V,G\right) $-closed submodule. Let $\left( Z_{0},H_{0}\right)
\leq \left( Z,H\right) $ be a finitely generated subrepresentation of the $%
\left( Z,H\right) $. $\left( Z_{0},H_{0}\right) \cong W\left( X,Y\right)
/\left( T,L\right) $, where $\left( T,L\right) $ is a normal
subrepresentation of $W\left( X,Y\right) $). By \textbf{Proposition 3.4} and
\textbf{Proposition 4.2}, $T$ is an action type $\left( Z,H\right) $-closed
submodule and an action type $\left( V,G\right) $-closed submodule. Hence,
by \textbf{Proposition 4.2} and by \textbf{Proposition 3.4}, there exists a
normal subrepresentation $\left( T,D\right) \leq W\left( X,Y\right) $ such
that $W\left( X,Y\right) /\left( T,D\right) \in \mathcal{SC}\left(
V,G\right) $. Therefore $\left( Z_{0},H_{0}\right) \in \mathcal{Q}^{0}%
\mathcal{Q}^{r}\mathcal{SC}\left( V,G\right) $ and $\left( Z,H\right) \in
\mathcal{LQ}^{0}\mathcal{Q}^{r}\mathcal{SC}\left( V,G\right) $. The proof is
complete.

We shall denote $\left( Z,H\right) \prec \left( V,G\right) $ if and only if $%
\left( Z,H\right) \in \mathcal{LQ}^{0}\mathcal{Q}^{r}\mathcal{SC}\left(
V,G\right) $. By consideration of action type closed submodules we have

\textbf{Proposition 5.1. }\textit{The relation "}$\prec $\textit{" is the
preorder in the }$Rep-K$\textit{.}

By \textbf{Corollary} from the \textbf{Proposition 4.1} we have

\textbf{Corollary 1} from the \textbf{Theorem 5.1.} \textit{Let }$\left(
V_{1},G_{1}\right) ,\left( V_{2},G_{2}\right) \in Rep-K$\textit{. }$\left(
V_{1},G_{1}\right) \sim _{a.t.}\left( V_{2},G_{2}\right) $\textit{\ if and
only if }$\left( V_{1},G_{1}\right) \prec \left( V_{2},G_{2}\right) $\textit{%
\ and }$\left( V_{2},G_{2}\right) \prec \left( V_{1},G_{1}\right) $\textit{\
i.e., if and only if }$\left( V_{1},G_{1}\right) \in \mathcal{LQ}^{0}%
\mathcal{Q}^{r}\mathcal{SC}\left( V_{2},G_{2}\right) $\textit{\ and }$\left(
V_{2},G_{2}\right) \in \mathcal{LQ}^{0}\mathcal{Q}^{r}\mathcal{SC}\left(
V_{1},G_{1}\right) $\textit{.}

By (2.11), (2.5) and monotony of operators: $\mathcal{Q}^{r}$, $\mathcal{Q}%
^{0}$, $\mathcal{L}$, $\mathcal{S}$, $\mathcal{C}$ we have

\textbf{Corollary 2} from the \textbf{Theorem 5.1.} \textit{Let }$\left(
V_{1},G_{1}\right) ,\left( V_{2},G_{2}\right) \in Rep-K$\textit{. }$\left(
V_{1},G_{1}\right) \sim _{a.t.}\left( V_{2},G_{2}\right) $\textit{\ if and
only if }$\mathcal{LQ}^{0}\mathcal{Q}^{r}\mathcal{SC}\left(
V_{1},G_{1}\right) =\mathcal{LQ}^{0}\mathcal{Q}^{r}\mathcal{SC}\left(
V_{2},G_{2}\right) $\textit{.}

\textbf{Corollary 3} from the \textbf{Theorem 5.1.} \textit{Let }$\left(
V,G\right) \in Rep-K$\textit{. Then}%
\begin{equation*}
\mathcal{LQ}^{0}\mathcal{Q}^{r}\mathcal{SC}\left( V,G\right) \subset
qVar_{a.t.}\left( V,G\right) .
\end{equation*}%
\textit{\ }

\textit{Proof:}

Let $\left( Z,H\right) \in \mathcal{LQ}^{0}\mathcal{Q}^{r}\mathcal{SC}\left(
V,G\right) $. If
\begin{equation*}
\left( V,G\right) \vDash \left( \left( \bigwedge\limits_{i\in
I}w_{i}=0\right) \Rightarrow \left( w_{0}=0\right) \right) ,
\end{equation*}%
then, by \textbf{Proposition 4.4}, $w_{0}\in \left\{ w_{i}\mid i\in
I\right\} _{\left( V,G\right) }^{\nabla \nabla }$\textit{. }Every action
type $\left( Z,H\right) $-closed submodule is also an action type $\left(
V,G\right) $-closed submodule, therefore $\left\{ w_{i}\mid i\in I\right\}
_{\left( Z,H\right) }^{\nabla \nabla }\ $$\supset \left\{ w_{i}\mid i\in
I\right\} _{\left( V,G\right) }^{\nabla \nabla }\ni w_{0}$. Thus,
\begin{equation*}
\left( Z,H\right) \vDash \left( \left( \bigwedge\limits_{i\in
I}w_{i}=0\right) \Rightarrow \left( w_{0}=0\right) \right) .
\end{equation*}%
The proof is complete.

\textbf{Corollary 4} from the \textbf{Theorem 5.1.} \textit{If two
representations }$\left( V_{1},G_{1}\right) $ ,$\left( V_{2},G_{2}\right)
\in Rep-K$\textit{\ are equivalent then they are action type geometrically
equivalent. In particular, every representation }$\left( V,G\right) \in
Rep-K $\textit{\ is action type geometrically equivalent to its faithful
image }$\left( V,\widetilde{G}\right) $\textit{.}

\textit{Proof:}

Let representations $\left( V_{1},G_{1}\right) $ and $\left(
V_{2},G_{2}\right) $ be equivalent; then
\begin{equation*}
\mathcal{Q}^{0}\mathcal{Q}^{r}\left( V_{1},G_{1}\right) =\mathcal{Q}^{0}%
\mathcal{Q}^{r}\left( V_{2},G_{2}\right) .
\end{equation*}%
So%
\begin{equation*}
\mathcal{LQ}^{0}\mathcal{Q}^{r}\mathcal{SC}\left( V_{1},G_{1}\right) \supset
\mathcal{Q}^{0}\mathcal{Q}^{r}\left( V_{2},G_{2}\right) \ni \left(
V_{2},G_{2}\right) ,
\end{equation*}%
i.e., $\left( V_{2},G_{2}\right) \prec \left( V_{1},G_{1}\right) $. By
symmetry the proof is complete.

\section{Action type Noetherianity of representations.}

\textbf{Definitions 6.1.}\textit{\ We call a variety }$\Theta \subset Rep-K$%
\textit{\ }\textbf{action type Noetherian}\textit{\ if for every finite }$X$%
\textit{\ \ and }$Y$\textit{\ every }$KF_{\Theta }\left( Y\right) $\textit{%
-submodule of }$E_{\Theta }(X,Y)$ is finitely generated (as a $KF_{\Theta
}\left( Y\right) $\textit{-module).}

It is clear that action type Noetherianity of the variety $\Theta $ is
equivalent to the ascending chain condition for $KF_{\Theta }\left( Y\right)
$-submodules of $E_{\Theta }\left( X,Y\right) $ for every $X$ and $Y$.
Action type geometrical Noetherianity of the representation $\left(
V,G\right) $ is equivalent to the ascending chain condition for action type $%
\left( V,G\right) $-closed submodules of $E_{\Theta }(X,Y)$. So, by \textbf{%
Proposition 4.3}, every representation $\left( V,G\right) $ from the action
type Noetherian variety $\Theta $ is action type geometrically Noetherian.

Similarly to [Pl4], Proposition 7 one can prove that:

\textbf{Proposition 6.1.} \textit{A representation }$\left( V,G\right) \in
\Theta $\textit{\ is action type logically Noetherian\ if and only if the
union of any directed system of action type }$\left( V,G\right) $\textit{%
-closed submodules of }$E_{\Theta }(X,Y)$\textit{\ for every }$X$\textit{\
and }$Y$\textit{\ is also an action type }$\left( V,G\right) $\textit{%
-closed submodule.}

Hence every action type geometrically Noetherian representation is also
action type logically Noetherian.

If $\left( V,G\right) $ is an action type logically Noetherian
representation, then, by \textbf{Proposition 4.4}, every infinite action
type quasi-identity (0.1') can be reduced to the finite action type
quasi-identity (0.1).

We shall give some examples of these notions.

\textbf{Proposition 6.2.} \textit{The variety }$\mathfrak{S}^{n}$\textit{\
of representation over the Noetherian ring }$K$\textit{\ is the action type
Noetherian variety for every }$n\in
\mathbb{N}
$\textit{.}

\textit{Proof: }

We denote $\bigoplus\limits_{x\in X}x\left( KF\left( Y\right) /\Delta
^{n}\right) =X\left( KF\left( Y\right) /\Delta ^{n}\right) $, where $\Delta $
is the augmentation ideal of the $KF\left( Y\right) $. $W_{\mathfrak{S}%
^{n}}\left( X,Y\right) \left( XKF\left( Y\right) /\bigoplus\limits_{x\in
X}x\Delta ^{n},F\left( Y\right) \right) \cong \left( X\left( KF\left(
Y\right) /\Delta ^{n}\right) ,F\left( Y\right) \right) $. If $\left\vert
Y\right\vert =m$, then, by the Taylor formula for Fox derivation ([Vvs]),%
\begin{equation*}
w=w^{\varepsilon }+\sum\limits_{i_{1}=1}^{m}\left( \partial _{i_{1}}w\right)
^{\varepsilon }\left( y_{i_{1}}-1\right)
+\sum\limits_{i_{1},i_{2}=1}^{m}\left( \partial _{i_{1}i_{2}}w\right)
^{\varepsilon }\left( y_{i_{1}}-1\right) \left( y_{i_{2}}-1\right) +\ldots
\end{equation*}%
\begin{equation*}
\ldots +\sum\limits_{i_{1},\ldots ,i_{k-1}=1}^{m}\left( \partial
_{i_{1},\ldots ,i_{k-1}}w\right) ^{\varepsilon }\left( y_{i_{1}}-1\right)
\ldots \left( y_{i_{k-1}}-1\right) +
\end{equation*}%
\begin{equation*}
\sum\limits_{i_{1},\ldots ,i_{k}=1}^{m}\left( \partial _{i_{1},\ldots
,i_{k}}w\right) \left( y_{i_{1}}-1\right) \ldots \left( y_{i_{k}}-1\right)
\end{equation*}%
for every $w\in KF\left( Y\right) $ and for every $k\in
\mathbb{N}
$, where $\varepsilon :KF\left( Y\right) \rightarrow K$ is the augmentation
homomorphism, $\partial _{i_{1},\ldots ,i_{s}}:$ $KF\left( Y\right)
\rightarrow KF\left( Y\right) $ is the $s$-th Fox derivation by the
variables $y_{i_{1}},\ldots ,y_{i_{s}}$ ($1\leq s\leq k$). So, $KF\left(
Y\right) /\Delta ^{n}$ is the finitely generated $K$-module for every $Y$.
Hence, $X\left( KF\left( Y\right) /\Delta ^{n}\right) $ is the finitely
generated $K$-module for every $X$ and $Y$. $K$ is the Noetherian ring, so,
every $K$-submodule and every $KF\left( Y\right) $-submodule of $X\left(
KF\left( Y\right) /\Delta ^{n}\right) $ is finitely generated. The proof is
complete.

\textbf{Proposition 6.3.} \textit{Every faithful finitely dimension
representation }$\left( V,G\right) $\textit{\ over the field }$K$\textit{\
is action type geometrically Noetherian.}

This proposition we can prove by using ideas from [BMR, Theorem B1].

\textbf{Corollary.} \textit{Every finite dimension representation }$\left(
V,G\right) $\textit{\ over the field }$K$\textit{\ is action type
geometrically Noetherian.}

\textit{Proof:} By \textbf{Corollary 4} from the \textbf{Theorem 5.1} and
\textbf{Corollary} from the \textbf{Proposition 4.1}.

\textbf{Theorem 6.1.} \textit{Let }$\left( V_{1},G_{1}\right) $\textit{\ and
}$\left( V_{2},G_{2}\right) $\textit{\ be action type logically Noetherian
representations. Then }$\left( V_{1},G_{1}\right) \sim _{a.t.}\left(
V_{2},G_{2}\right) $\textit{\ if and only if }$qId_{a.t.}\left(
V_{1},G_{1}\right) =qId_{a.t.}\left( V_{2},G_{2}\right) $\textit{.}

\textit{Proof:} By \textbf{Proposition 4.4}, \textbf{Corollary 1} from
\textbf{Proposition 4.4} and \textbf{Corollary} from \textbf{Proposition 4.1}%
.

\textbf{Corollary.} \textit{In an action type Noetherian variety of
representations }$\Theta $\textit{\ there is bijection between classes of
action type geometrical equivalent representations and action type
quasi-varieties generated by one representation.}

\textit{Proof:}

Let $\left( V,G\right) \in \Theta $. We denote by $\left[ \left( V,G\right) %
\right] $ the class of all representations in $\Theta $ which are action
type geometrically equivalent to the representation $\left( V,G\right) $. It
is easy to check that the correspondence $\left[ \left( V,G\right) \right]
^{\varphi }=qVar_{a.t.}\left( V,G\right) $ is well defined, and bijection.

\textbf{Proposition 6.4.} \textit{Let }$\left( V,G\right) $ \textit{be an
action type logically Noetherian representation. Then }$qVar_{a.t.}\left(
V,G\right) \subset \mathcal{LQ}^{0}\mathcal{Q}^{r}\mathcal{SC}\left(
V,G\right) $\textit{.}

\textit{Proof:}

Let $\left( Z,H\right) \in qVar_{a.t.}\left( V,G\right) $ and $T\leq \left(
XKF\left( Y\right) \right) _{KF\left( Y\right) }$ is the action type $\left(
Z,H\right) $-closed submodule, but not action type $\left( V,G\right) $%
-closed submodule. Let $w\in T_{\left( V,G\right) }^{\nabla \nabla
}\backslash T$. There is $T_{0}=\left\{ w_{1},\ldots ,w_{n}\right\} \subset
T $, such that $w\in \left( T_{0}\right) _{\left( V,G\right) }^{\nabla
\nabla } $. Therefore,%
\begin{equation*}
\left( V,G\right) \vDash \left( \left( \bigwedge\limits_{i=1}^{n}\left(
w_{i}=0\right) \right) \Rightarrow \left( w=0\right) \right)
\end{equation*}%
and%
\begin{equation*}
\left( Z,H\right) \vDash \left( \left( \bigwedge\limits_{i=1}^{n}\left(
w_{i}=0\right) \right) \Rightarrow \left( w=0\right) \right) .
\end{equation*}%
So $w\in \left( T_{0}\right) _{\left( Z,H\right) }^{\nabla \nabla }\subset
T_{\left( Z,H\right) }^{\nabla \nabla }$, but $w\notin T$. By this
contradiction, $T$ is the action type $\left( V,G\right) $-closed submodule.
By \textbf{Theorem 5.1}, the proof is complete.

\textbf{Theorem 6.2.} $\mathcal{LQ}^{0}\mathcal{Q}^{r}$\QTR{cal}{SC}$\left(
V,G\right) =qVar_{a.t.}\left( V,G\right) $\textit{\ if and only if }$\left(
V,G\right) $\textit{\ is an action type logically Noetherian representation.}

\textit{Proof: }

By \textbf{Corollary 3} from the \textbf{Theorem 5.1}, we always have
\begin{equation*}
\mathcal{LQ}^{0}\mathcal{Q}^{r}\mathcal{SC}\left( V,G\right) \subset
qVar_{a.t.}\left( V,G\right) .
\end{equation*}%
If $\left( V,G\right) $ is an action type logically Noetherian
representation, then, by \textbf{Proposition 6.4},
\begin{equation*}
\mathcal{LQ}^{0}\mathcal{Q}^{r}\mathcal{SC}\left( V,G\right)
=qVar_{a.t.}\left( V,G\right) .
\end{equation*}

Let%
\begin{equation*}
\mathcal{LQ}^{0}\mathcal{Q}^{r}\mathcal{SC}\left( V,G\right)
=qVar_{a.t.}\left( V,G\right) .
\end{equation*}%
Let $\left\{ T_{i}\mid i\in I\right\} $ be a direct system of action type $%
\left( V,G\right) $-closed submodules of $XKF\left( Y\right) $ and $%
T=\bigcup\limits_{i\in I}T_{i}$.

Let $\left( V,G\right) \vDash \mathfrak{q}$ where $\mathfrak{q}$ is an
action type quasi-identity. By \textbf{Proposition 4.2} and \textbf{%
Proposition 3.4}%
\begin{equation*}
\left( XKF\left( Y\right) /T_{i},F\left( Y\right) \right) \in \mathcal{Q}^{0}%
\mathcal{SC}\left( V,G\right) \subset qVar_{a.t.}\left( V,G\right)
\end{equation*}%
for every $i\in I$, so, using the method of [Pl4, Theorem 1], we can prove
that $\left( XKF\left( Y\right) /T,F\left( Y\right) \right) \vDash $ $%
\mathfrak{q}$. Hence
\begin{equation*}
\left( XKF\left( Y\right) /T,F\left( Y\right) \right) \in qVar_{a.t.}\left(
V,G\right) =\mathcal{LQ}^{0}\mathcal{Q}^{r}\mathcal{SC}\left( V,G\right) .
\end{equation*}%
Consequently, there exists $H\trianglelefteq F\left( Y\right) $, such that
\begin{equation*}
\left( XKF\left( Y\right) /T,F\left( Y\right) /H\right) \in \mathcal{SC}%
\left( V,G\right) .
\end{equation*}%
Therefore, by \textbf{Proposition 3.4} and \textbf{Proposition 4.2} $T$ is
an action type $\left( V,G\right) $-closed submodule. The proof is complete.

\section{Action type quasi-varieties of representations.}

\textbf{Definition 7.1.}\textit{\ We say that a class of representations }$%
\mathfrak{X}$\textit{\ is right hereditary if }$\mathcal{S}_{r}\mathfrak{X}=%
\mathfrak{X}$\textit{.}

In [Ma] it was proved that a class of algebras which contains the unit
algebra is a quasi-variety if and only if this class is closed under the
operators $\mathcal{S}$ and $\mathcal{F}$. Later on in [Gv] this result was
established for the case of many sorted algebras. We use this fact in order
to describe the action type quasi-varieties\ of representations. It is clear
that every non empty class of representations which is closed under the
operators $\mathcal{S}$ contains the unit representation $\left( \left\{
0\right\} ,\left\{ 1\right\} \right) $, so the non empty class of
representations is a quasi-variety if and only if this class is closed under
the operators $\mathcal{S}$ and $\mathcal{F}$.

Let $\mathfrak{X}$ be a class of representations. Denote by $\mathfrak{X}%
_{G} $ the class of all $KG$-modules $V_{KG}$, such that the corresponding
representation $\left( V,G\right) $ belongs to the class $\mathfrak{X}$.

\textbf{Lemma 7.1.} A class \textit{\ }$\mathfrak{X}$\textit{\ is a
saturated quasi-variety of representations if and only if }$\mathfrak{X}$%
\textit{\ is saturated, right hereditary and }$\mathfrak{X}_{G}$\textit{\ is
a quasi-variety of }$KG$\textit{-modules for every group }$G$\textit{.}

\textit{Proof: }

Let $\mathfrak{X}$ be a saturated quasi-variety. It is clear that $\mathfrak{%
X}$ is a right hereditary class.

If $G$ is a group and $W$ is a submodule of the $KG$-module $V$, then $%
\left( W,G\right) $ is a subrepresentation of the $\left( V,G\right) $. So $%
\mathcal{S}\mathfrak{X}_{G}\subset \mathfrak{X}_{G}$.

Let $\left\{ \left( V_{i}\right) _{KG}\mid i\in I\right\} \subset \mathfrak{X%
}_{G}$, $\mathfrak{F}$ be a filter in the $I$. $\left( V_{i},G\right) \in
\mathfrak{X}$ for every $i\in I$. The filtered product of the family $%
\left\{ \left( V_{i}\right) _{KG}\mid i\in I\right\} $ as $KG$-modules is $%
\left( \dprod\limits_{i\in I}\left( V_{i}\right) _{KG}\right) /\sim _{%
\mathfrak{F}}$. The filtered product of representations $\left\{ \left(
V_{i},G\right) \mid i\in I\right\} $ is $\left( \left( \dprod\limits_{i\in
I}V_{i}\right) /\sim _{\mathfrak{F}},\left( G^{I}\right) /\sim _{\mathfrak{F}%
}\right) $. The representation $\left( \left( \dprod\limits_{i\in
I}V_{i}\right) /\sim _{\mathfrak{F}},G\right) $ is its subrepresentation,
because the diagonal of $\left( G^{I}\right) /\sim _{\mathfrak{F}}$is
isomorphic to $G$, so $\left( \dprod\limits_{i\in I}\left( V_{i}\right)
_{KG}\right) /\sim _{\mathfrak{F}}\in \mathfrak{X}_{G}$. Therefore $\mathcal{%
F}\mathfrak{X}_{G}\subset \mathfrak{X}_{G}$ and, by [Ma], $\mathfrak{X}_{G}$
is a quasi-variety of $KG$-modules for every group $G$.

Let $\mathfrak{X}$ be saturated, right hereditary and $\mathfrak{X}_{G}$ be
a quasi-variety of $KG$-modules for every group $G$. Let $\left( W,H\right)
\leq \left( V,G\right) $, $\left( V,G\right) \in \mathfrak{X}$. Then $%
V_{KH}\in \mathfrak{X}_{H}$ and $W_{KH}\in \mathcal{S}\mathfrak{X}_{H}=%
\mathfrak{X}_{H}$. Therefore $\mathcal{S}\mathfrak{X}=\mathfrak{X}$.

Let $\left\{ \left( V_{i},G_{i}\right) \mid i\in I\right\} \subset \mathfrak{%
X}$, $\mathfrak{F}$ be a filter over the $I$. Denote $\dprod\limits_{i\in
I}G_{i}=G$. The filtered product of the family of representations $\left\{
\left( V_{i},G_{i}\right) \mid i\in I\right\} $ is $\left( \left(
\dprod\limits_{i\in I}V_{i}\right) /\sim _{\mathfrak{F}},G/\sim _{\mathfrak{F%
}}\right) $. Let $\pi _{i}:G\rightarrow G_{i}$ be projections. Epimorphism $%
\pi _{i}$ defines representation $\left( V_{i},G\right) $ and $\left(
id_{V_{i}},\pi _{i}\right) $ is a homomorphism of representations for every $%
i\in I$. Hence, $\left( V_{i},G\right) \in \mathcal{Q}^{0}\mathfrak{X\subset
X}$ and $\left( V_{i}\right) _{KG}\in \mathfrak{X}_{G}$ for every $i\in I$.
So $\left( \dprod\limits_{i\in I}\left( V_{i}\right) _{KG}\right) /\sim _{%
\mathfrak{F}}\in \mathfrak{X}_{G}$ and $\left( \left( \dprod\limits_{i\in
I}V_{i}\right) /\sim _{\mathfrak{F}},G\right) \in \mathfrak{X}$. So, $\left(
\left( \dprod\limits_{i\in I}V_{i}\right) /\sim _{\mathfrak{F}},G/\sim _{%
\mathfrak{F}}\right) \in \mathcal{Q}^{r}\mathfrak{X\subset X}$. Therefore $%
\mathcal{F}\mathfrak{X}\subset \mathfrak{X}$ and, by [Gv], $\mathfrak{X}$ is
a quasi-variety. The proof is complete.

\textbf{Theorem 7.1.}\textit{\ A quasi-variety of representations }$%
\mathfrak{X}$\textit{\ is an action type quasi-variety of representations if
and only if }$\mathfrak{X}$\textit{\ is a saturated quasi-variety.}

\textit{Proof:}

Let $\mathfrak{X}$ be an action type quasi-variety of representations. Let $%
\left( V,G\right) \in \mathfrak{X}$ and a representation $\left( W,H\right) $
is equivalent to the representation $\left( V,G\right) $. By \textbf{%
Corollary 4} from the \textbf{Theorem 5.1}, $\left( V,G\right) \sim
_{a.t.}\left( W,H\right) $ and, by \textbf{Corollary 1} from \textbf{%
Proposition 4.4}, $qId_{a.t.}\left( V,G\right) =qId_{a.t.}\left( W,H\right) $%
, so $\left( W,H\right) \in \mathfrak{X}$. Therefore $\mathfrak{X}$ is a
saturated class of representations.

Let $\mathfrak{X}$\textit{\ }be a saturated quasi-variety of
representations. By \textbf{Lemma 7.1}, $\mathfrak{X}_{F_{\infty }}$ is a
quasi-variety of $KF_{\infty }$-modules, i.e., $\mathfrak{X}_{F_{\infty
}}=qVar\mathfrak{Q}$, where $\mathfrak{Q}=\left\{ \mathfrak{q}_{i}\mid i\in
I\right\} ,$%
\begin{equation}
\mathfrak{q}_{i}\equiv \left( \forall x_{1}\ldots \forall x_{n_{i}}\left(
\left( \bigwedge\limits_{j=1}^{k_{i}}\left( w_{ij}=0\right) \right)
\Rightarrow \left( w_{i0}=0\right) \right) \right) ,  \tag{7.1}
\end{equation}%
$w_{ij}=w_{ij}\left( x_{1},\ldots ,x_{n_{i}},y_{1},\ldots ,y_{m_{i}}\right)
\in
\mathbb{N}
KF_{\infty }$ ($F_{\infty }$ is a free group with the countable set of
generators $\left\{ y_{1},y_{2},\ldots \right\} $, $%
\mathbb{N}
KF_{\infty }$ is the free $KF_{\infty }$-module with the countable basis $%
\left\{ x_{1},x_{2},\ldots \right\} $). In (7.1) we consider $%
x_{1},x_{2},\ldots $ as variables and $y_{1},y_{2},\ldots $ as constants.
But we can consider $y_{1},y_{2},\ldots $ also as variables. By this point
of view, $\widetilde{\mathfrak{q}_{i}}=\forall y_{1}\ldots \forall y_{m_{i}}%
\mathfrak{q}_{i}$ is an action type quasi-identity in $Rep-K$ and the set $%
\widetilde{\mathfrak{Q}}=\left\{ \widetilde{\mathfrak{q}_{i}}\mid i\in
I\right\} $ will be a set of action type quasi-identities in $Rep-K$. We
shall prove that $\mathfrak{X}=qVar\widetilde{\mathfrak{Q}}$.

Let $\left( V,G\right) \vDash \widetilde{\mathfrak{Q}}$. Let $G_{0}\leq G$
be a finitely generated subgroup of the group $G$. Also, $\left(
V,G_{0}\right) \vDash \widetilde{\mathfrak{Q}}$. There is an epimorphism $%
\beta :F_{\infty }\rightarrow G_{0}$. Denote by $\left( V\right) _{\beta }$
the $KF_{\infty }$-module defined by the homomorphism $\beta $ and by $%
\left( V,F_{\infty }\right) $ the representation corresponding to the module
$\left( V\right) _{\beta }$. Let $\widetilde{\mathfrak{q}_{i}}\in \widetilde{%
\mathfrak{Q}}$. The mapping $\alpha :\left\{ x_{l}\mid l\in
\mathbb{N}
\right\} \rightarrow V$ can be extended to the homomorphism of $KF_{\infty }$%
-modules $\alpha :%
\mathbb{N}
KF_{\infty }\rightarrow \left( V\right) _{\beta }$. It is clear that in this
situation the pair $\left( \alpha ,\beta \right) $ will be a homomorphism of
representations: $\left( \alpha ,\beta \right) :\left(
\mathbb{N}
KF_{\infty },F_{\infty }\right) \rightarrow \left( V,G_{0}\right) $. The
result $w_{ij}^{\alpha }$ does not depend on point of view on $\alpha $: as
a homomorphism of $KF_{\infty }$-modules or as a left component of
homomorphism of representations. So $\left( V,G_{0}\right) \vDash \widetilde{%
\mathfrak{q}_{i}}$ if and only if $\left( V\right) _{\beta }\vDash \mathfrak{%
q}_{i}$. Hence, $\left( V\right) _{\beta }\vDash \mathfrak{Q}$. Therefore, $%
\left( V\right) _{\beta }\in \mathfrak{X}_{F_{\infty }}$ and $\left(
V,G_{0}\right) \in \mathcal{Q}^{r}\mathfrak{X\subset }$ $\mathfrak{X}$. $%
\mathfrak{X}$ is a quasi-variety, so $\left( V,G\right) \in \mathfrak{X}$,
because all quasi-identities which define $\mathfrak{X}$ are checked in
finitely generated representations.

Let $\left( V,G\right) \in \mathfrak{X}$. Let $\left( \alpha ,\beta \right)
:\left(
\mathbb{N}
KF_{\infty },F_{\infty }\right) \rightarrow \left( V,G\right) $ be a
homomorphism of representation. Denote by $\left( V\right) _{\beta }$ the $%
KF_{\infty }$-module, defined by the homomorphism $\beta $, and $\left(
V,F_{\infty }\right) $ - the representation corresponding to this module. We
have that $\left( V,\mathrm{im}\beta \right) \in \mathfrak{X}$, $\left(
V,F_{\infty }\right) \in \mathcal{Q}^{0}\mathfrak{X\subset }$ $\mathfrak{X}$
and $\left( V\right) _{\beta }\in \mathfrak{X}_{F_{\infty }}$. Therefore $%
\left( V\right) _{\beta }\vDash \mathfrak{Q}$. Because $\alpha :%
\mathbb{N}
KF_{\infty }\rightarrow \left( V\right) _{\beta }$ is a homomorphism of $%
KF_{\infty }$-modules, as above, $\left( V,G\right) \vDash \widetilde{%
\mathfrak{Q}}$. The proof is complete.

\textbf{Corollary 1.} $qVar_{a.t.}\mathfrak{X}=\mathcal{Q}^{0}\mathcal{Q}^{r}%
\mathcal{SF}\mathfrak{X}$\textit{\ for every class of representations }$%
\mathfrak{X}$\textit{.}

\textit{Proof:}

$qVar_{a.t.}\mathfrak{X}=qVar\left( qId_{a.t.}\mathfrak{X}\right) $ is a
quasi-variety, so $qVar_{a.t.}\mathfrak{X}$ is closed by $\mathcal{S}$ and $%
\mathcal{F}$. By \textbf{Theorem 7.1}, $qVar_{a.t.}\mathfrak{X}$ is a
saturated class of representations. Thus it is closed under $\mathcal{Q}^{0}$
and $\mathcal{Q}^{r}$. $\mathfrak{X}\subset qVar_{a.t.}\mathfrak{X}$ thus,
by (2.12),%
\begin{equation*}
\mathcal{Q}^{0}\mathcal{Q}^{r}\mathcal{SF}\mathfrak{X}\subset qVar_{a.t.}%
\mathfrak{X}.
\end{equation*}

Also, by (2.12), $\mathcal{Q}^{0}\mathcal{Q}^{r}\mathcal{SF}\mathfrak{X}$ is
closed under $\mathcal{S}$, $\mathcal{F}$, $\mathcal{Q}^{0}$ and $\mathcal{Q}%
^{r}$. By \textbf{Theorem 7.1}, $\mathcal{Q}^{0}\mathcal{Q}^{r}\mathcal{SF}%
\mathfrak{X}$ is\ an action type quasi-variety of representations, i.e.,%
\begin{equation*}
\mathcal{Q}^{0}\mathcal{Q}^{r}\mathcal{SF}\mathfrak{X}=qVar\mathfrak{Q},
\end{equation*}%
where $\mathfrak{Q}$ is the set of action type quasi-identities. $\mathcal{Q}%
^{0}$, $\mathcal{Q}^{r}$, $\mathcal{S}$, $\mathcal{F}$ are operators of
extension, so%
\begin{equation*}
\mathfrak{X}\subset \mathcal{Q}^{0}\mathcal{Q}^{r}\mathcal{SF}\mathfrak{X}%
=qVar\mathfrak{Q}.
\end{equation*}%
Hence $\mathfrak{X}\vDash \mathfrak{Q}$ and $qId_{a.t.}\mathfrak{X}\supset
\mathfrak{Q}$. Therefore%
\begin{equation*}
qVar_{a.t.}\mathfrak{X}\subset qVar\mathfrak{Q}=\mathcal{Q}^{0}\mathcal{Q}%
^{r}\mathcal{SF}\mathfrak{X}.
\end{equation*}%
The proof is complete.

\textbf{Corollary 2.} $qVar_{a.t.}\mathfrak{X}=\mathcal{Q}^{0}\mathcal{Q}^{r}%
\mathcal{SCC}_{up}\mathfrak{X}$\textit{\ for every class of representations }%
$\mathfrak{X}$\textit{.}

\textit{Proof:} By [GL], for every class of algebras $\mathfrak{X}$ we have $%
\mathcal{F}\mathfrak{X}=\mathcal{CC}_{up}\mathfrak{X}$.

\section{Existing of continuum non isomorphic simple $F_{2}$-modules.}

There is a continuum of non isomorphic simple $2$-generated groups ([Ca]).
Using this fact, R.Gobel and S. Shelah proved [GSh] that there is a non
logically Noetherian group.

$K$ in this and the next section is a field such that $\left\vert
K\right\vert =\aleph _{0}$ and $charK\neq 2$. We shall prove that there is
continuum of non isomorphic simple modules over $KF_{2}$, where $F_{2}$ is
the free group with $2$ generators. And we shall deduce from this fact that
there is a non action type logically Noetherian representation.

Let $\Delta $ be the augmentation ideal of group algebra $KG$.

\textbf{Proposition 8.1.} \textit{\ Let }$G$\textit{\ be a non periodic
group, i.e. }$\exists $\textit{\ }$g\in G$\textit{, such that }$\left\vert
g\right\vert =\infty $\textit{. Then the set }$\Omega =$\textit{\ }$\left\{
U\lneqq KG_{KG}\mid U\nsubseteq \Delta \right\} $\textit{\ is non empty and
has a maximal element }$U_{G}$\textit{\ which is a maximal right ideal in }$%
KG$\textit{.}

\textit{Proof: \ }

The element $g+1$ is not invertible in $KG$, because $\left\vert
g\right\vert =\infty $. So, $\left( g+1\right) KG\lneqq KG_{KG}$. Also $%
\left( g+1\right) KG\nsubseteq \Delta $. So, $\Omega \neq \varnothing $. By
Zorn's lemma the set $\Omega $ has a maximal element $U_{G}$. It is easy to
check that $U_{G}$ is a maximal right ideal in $KG$.

\textbf{Corollary.} $KG/U_{G}$\textit{\ is a simple }$KG$\textit{-module in
the situation of }\textbf{Proposition 8.1.}

The ideal $U_{G}$, of course, is not uniquely defined by the group $G$.

\textbf{Proposition 8.2.} \textit{If }$\Gamma $\textit{\ is a simple non
periodic group, then the representation }$\left( K\Gamma /\mathrm{ann}%
_{K\Gamma }\left( K\Gamma /U_{\Gamma }\right) ,\Gamma \right) $\textit{\ is
faithful.}

\textit{Proof: }

$\mathrm{ann}_{K\Gamma }\left( K\Gamma /U_{\Gamma }\right) =L_{\Gamma
}\subset U_{\Gamma }$, by \textbf{Proposition 1.1}. We consider two
representations: $\left( K\Gamma /L_{\Gamma },\Gamma \right) $ and $\left(
K\Gamma /U_{\Gamma },\Gamma \right) $.%
\begin{equation*}
\left( 1+U_{\Gamma }\right) \cap \Gamma =\left\{ g\in \Gamma \mid g-1\in
U_{\Gamma }\right\} \supset \left( 1+L_{\Gamma }\right) \cap \Gamma .
\end{equation*}%
$\left\{ g\in \Gamma \mid g-1\in U_{\Gamma }\right\} \neq \Gamma $,
otherwise $U_{\Gamma }=\Delta $. By \textbf{Corollary} from \textbf{%
Proposition 1.2},%
\begin{equation*}
\left( 1+L_{\Gamma }\right) \cap \Gamma =\ker \left( K\Gamma /L_{\Gamma
},\Gamma \right) \vartriangleleft \Gamma
\end{equation*}%
$\Gamma $ is a simple group, so $\ker \left( K\Gamma /L_{\Gamma },\Gamma
\right) =\left\{ 1\right\} $. The proof is complete.

\textbf{Corollary.} \textit{The group }$\Gamma $\textit{\ is embedded into
the associative algebra }\newline
$K\Gamma /\mathrm{ann}_{K\Gamma }\left( K\Gamma /U_{\Gamma }\right) $
\textit{in the situation of \textbf{Proposition 8.2}.}

\textit{Proof:}%
\begin{equation*}
\ker \left( K\Gamma /\mathrm{ann}_{K\Gamma }\left( K\Gamma /U_{\Gamma
}\right) ,\Gamma \right) =\left\{ g\in \Gamma \mid g-1\in \mathrm{ann}%
_{K\Gamma }\left( K\Gamma /U_{\Gamma }\right) \right\} =\left\{ 1\right\} .
\end{equation*}

\textbf{Theorem 8.1.} \textit{There exists a continuum of non isomorphic
simple modules over }$KF_{2}$\textit{, where }$F_{2}$\textit{\ is a free
group with }$2$\textit{\ generators.}

\textit{Proof: }

Let $\Re $ be the set of all non isomorphic simple $2$-generated groups,
considered in [Ca] ($\left\vert \Re \right\vert =\aleph $). By constructions
of [Ca], every $\Gamma \in \Re $ is a non periodic group. So, by \textbf{%
Proposition 8.1}, we can choose for every $\Gamma \in \Re $ the maximal
right ideal in $K\Gamma $: $U_{\Gamma }<K\Gamma _{K\Gamma }$. It holds that $%
U_{\Gamma }\neq \Delta $. $K\Gamma /U_{\Gamma }$ is the simple $K\Gamma $%
-module and the simple $KF_{2}$-module defined by the natural homomorphism $%
F_{2}\rightarrow \Gamma $. After the choosing of $U_{\Gamma }$\ for every $%
\Gamma \in \Re $, we define in $\Re $ the equivalence: $\Gamma _{1}\approx
\Gamma _{2}$ if $K\Gamma _{1}/U_{\Gamma _{1}}\cong K\Gamma _{2}/U_{\Gamma
_{2}}$ as $KF_{2}$-modules ($\Gamma _{1},\Gamma _{2}\in \Re $).

If $K\Gamma _{1}/U_{\Gamma _{1}}\cong K\Gamma _{2}/U_{\Gamma _{2}}$ as $%
KF_{2}$-modules,\ then, by \textbf{Proposition 1.4},%
\begin{equation*}
\mathrm{ann}_{KF_{2}}\left( K\Gamma _{1}/U_{\Gamma _{1}}\right) =\mathrm{ann}%
_{KF_{2}}\left( K\Gamma _{2}/U_{\Gamma _{2}}\right) ,
\end{equation*}%
and by \textbf{Corollary 2} from \textbf{Proposition 1.3},%
\begin{equation*}
K\Gamma _{1}/\mathrm{ann}_{K\Gamma _{1}}\left( K\Gamma _{1}/U_{\Gamma
_{1}}\right) \cong K\Gamma _{2}/\mathrm{ann}_{K\Gamma _{2}}\left( K\Gamma
_{2}/U_{\Gamma _{2}}\right)
\end{equation*}%
as associative algebras. By the \textbf{Corollary} from \textbf{Proposition
8.2}%
\begin{equation*}
\Gamma _{1}\hookrightarrow K\Gamma _{1}/\mathrm{ann}_{K\Gamma _{1}}\left(
K\Gamma _{1}/U_{\Gamma _{1}}\right) \cong K\Gamma _{2}/\mathrm{ann}_{K\Gamma
_{2}}\left( K\Gamma _{2}/U_{\Gamma _{2}}\right) ,
\end{equation*}%
so $\Gamma _{1}$ is isomorphic to one of the multiplicative subgroup of the
associative algebra $K\Gamma _{2}/\mathrm{ann}_{K\Gamma _{2}}\left( K\Gamma
_{2}/U_{\Gamma _{2}}\right) $. $\left\vert K\Gamma _{2}\right\vert =\aleph
_{0}$, so $\left\vert K\Gamma _{2}/\mathrm{ann}_{K\Gamma _{2}}\left( K\Gamma
_{2}/U_{\Gamma _{2}}\right) \right\vert \leq \aleph _{0}$ and there is a
countable set of $2$-generated subgroups of \newline
$K\Gamma _{2}/\mathrm{ann}_{K\Gamma _{2}}\left( K\Gamma _{2}/U_{\Gamma
_{2}}\right) $. Therefore, the cardinality of classes by equivalence "$%
\approx $" is not bigger than $\aleph _{0}$. So, there are $\aleph $ classes
by equivalence "$\approx $". The proof is complete.

\section{Non action type logically Noetherian representation of the group $%
F_{2}$.}

In this section we shall prove that there is a non action type logically
Noetherian representation. Let $\mathcal{P}\subseteq \Re $ be the set of all
non isomorphic simple $2$-generated groups such that simple $KF_{2}$-modules
$\left\{ K\Gamma /U_{\Gamma }\mid \Gamma \in \mathcal{P}\right\} $ are non
isomorphic. By the \textbf{Theorem 8.1}, $\left\vert \mathcal{P}\right\vert
=\aleph $.

If $\varphi _{\Gamma }:F_{2}\rightarrow F_{2}/H=\Gamma \in \mathcal{P}$ is
the natural homomorphism of groups, then, by \textbf{Proposition 1.3}, $%
K\Gamma /U_{\Gamma }\cong KF_{2}/U_{\Gamma }^{\varphi _{\Gamma }^{-1}}$.
Denote $KF_{2}/U_{\Gamma }^{\varphi _{\Gamma }^{-1}}=V_{\Gamma }$. $%
V_{\Gamma }$ is a simple $KF_{2}$-module.

Let $\left\{ V_{j}\mid j\in J\right\} $ be the set of all finitely generated
right ideals in $KF_{2}$. $V=\prod\limits_{j\in J}\left( KF_{2}/V_{j}\right)
$ is the $KF_{2}$-module. So, we can consider the representation $\left(
V,F_{2}\right) $. $\left\vert KF_{2}/V_{j}\right\vert =\aleph _{0}$ for
every $j\in J$, $\left\vert J\right\vert =\aleph _{0}$, so $\left\vert
V\right\vert =\aleph _{0}$.

\textbf{Theorem 9.1.} \textit{The representation }$\left( V,F_{2}\right) $%
\textit{\ is non action type logically Noetherian.}

\textit{Proof: }

We shall prove that there is $\Gamma \in \mathcal{P}$ such that $U_{\Gamma
}^{\varphi _{\Gamma }^{-1}}<\left( KF_{2}\right) _{KF_{2}}$ is not action
type $\left( V,F_{2}\right) $-closed. Let $\Gamma \in \mathcal{P}$ and $%
U_{\Gamma }^{\varphi _{\Gamma }^{-1}}$ be the action type $\left(
V,F_{2}\right) $-closed. By \textbf{Proposition 4.2} and \textbf{Proposition
3.4}, there exists $H\trianglelefteq F_{2}$ such that $\left( U_{\Gamma
}^{\varphi _{\Gamma }^{-1}},H\right) \trianglelefteq \left(
KF_{2},F_{2}\right) $ and $\left( V_{\Gamma },F_{2}/H\right) \in \mathcal{SC}%
\left( V,F_{2}\right) $. So, there exists a homomorphism of representations $%
\left( \iota ,\eta \right) :\left( V_{\Gamma },F_{2}\right) \rightarrow
\left( V^{I},F_{2}^{I}\right) $ ($I$ is the set of indices), such that $%
\iota $ is a monomorphism. Since $V_{\Gamma }$ is a simple $KF_{2}$-module,
we can conclude, that there exists an embedding of $KF_{2}$-module $%
V_{\Gamma }\hookrightarrow \left( V\right) _{\widetilde{\eta }}$, where $%
\widetilde{\eta }$ is an endomorphism of $F_{2}$. $\left\vert V\right\vert
=\aleph _{0}$, $\left\vert \mathrm{End}\left( F_{2}\right) \right\vert
=\aleph _{0}$ (every endomorphism is defined by values on generators). In
the module $\left( V\right) _{\widetilde{\eta }}$ there is a countable set
of simple submodules (every simple submodule is a cyclic, so it is defined
by a generator). So only the countable set of modules $V_{\Gamma }$ can be
embedded into the modules of the kind $\left( V\right) _{\widetilde{\eta }}$%
. Therefore, by \textbf{Theorem 8.1}, there is $\Gamma _{0}\in \mathcal{P}$
such that the right side ideal $U_{\Gamma _{0}}^{\varphi _{\Gamma
_{0}}^{-1}} $ is not action type $\left( V,F_{2}\right) $-closed.

On the other hand, by \textbf{Proposition 4.2} and \textbf{Proposition 3.4},
$V_{j}$ is action type $\left( V,F_{2}\right) $-closed for every $j\in J$.
Therefore,%
\begin{equation*}
\left\{ V_{j}\mid j\in J,V_{j}\subseteq U_{\Gamma _{0}}^{\varphi _{\Gamma
_{0}}^{-1}}\right\} =\left\{ V_{j}\mid j\in J_{0}\right\}
\end{equation*}%
is the direct system of action type $\left( V,F_{2}\right) $-closed modules,
which unit $\bigcup\limits_{j\in J_{0}}V_{j}=U_{\Gamma _{0}}^{\varphi
_{\Gamma _{0}}^{-1}}$ is not a action type $\left( V,F_{2}\right) $-closed
module. So, the representation $\left( V,F_{2}\right) $ is non action type
logically Noetherian. The proof is complete.

\textbf{Corollary.} \textit{There exists }$\left( \widetilde{V},\widetilde{%
F_{2}}\right) $\textit{\ an ultrapower of }$\left( V,F_{2}\right) $\textit{\
which is not action type geometrically equivalent to the }$\left(
V,F_{2}\right) $\textit{.}

\textit{Proof:}

If $\mathcal{C}_{up}\left( V,F_{2}\right) \subset \mathcal{LQ}^{0}\mathcal{Q}%
^{r}\mathcal{SC}\left( V,F_{2}\right) $ then, by \textbf{Corollary 2} from
\textbf{Theorem 7.1}\ and (2.11),%
\begin{equation*}
qVar_{a.t.}\left( V,F_{2}\right) =\mathcal{Q}^{0}\mathcal{Q}^{r}\mathcal{SCC}%
_{up}\left( V,F_{2}\right) \subset \mathcal{LQ}^{0}\mathcal{Q}^{r}\mathcal{SC%
}\left( V,F_{2}\right) .
\end{equation*}%
So, by the \textbf{Corollary 3} from \textbf{Theorem 5.1},
\begin{equation*}
qVar_{a.t.}\left( V,F_{2}\right) =\mathcal{LQ}^{0}\mathcal{Q}^{r}\mathcal{SC}%
\left( V,F_{2}\right)
\end{equation*}%
and, by \textbf{Theorem 6.2}, $\left( V,F_{2}\right) $ is action type
logically Notherian. By this contradiction, there exists $\left( \widetilde{V%
},\widetilde{F_{2}}\right) $ an ultrapower of $\left( V,F_{2}\right) $, such
that%
\begin{equation*}
\left( \widetilde{V},\widetilde{F_{2}}\right) \notin \mathcal{LQ}^{0}%
\mathcal{Q}^{r}\mathcal{SC}\left( V,F_{2}\right) .
\end{equation*}
On the other hand,
\begin{equation*}
\left( \widetilde{V},\widetilde{F_{2}}\right) \in \mathcal{LQ}^{0}\mathcal{Q}%
^{r}\mathcal{SC}\left( \widetilde{V},\widetilde{F_{2}}\right) ,
\end{equation*}%
so%
\begin{equation*}
\mathcal{LQ}^{0}\mathcal{Q}^{r}\mathcal{SC}\left( \widetilde{V},\widetilde{%
F_{2}}\right) \neq \mathcal{LQ}^{0}\mathcal{Q}^{r}\mathcal{SC}\left(
V,F_{2}\right) ,
\end{equation*}%
and, by \textbf{Corollary 2} from \textbf{Theorem 5.1}, $\left( \widetilde{V}%
,\widetilde{F_{2}}\right) \nsim _{a.t.}\left( V,F_{2}\right) $. The proof is
complete.

\begin{center}
{\Large References:}
\end{center}

[Bi] Birkhoff G. \textit{On the structure of abstract algebras,} Proc.
Cambr. Phil. Soc. 31 (1935), 433-454.

[BMR] G. Baumslag, A. Myasnikov, V. Remeslennikov. \textit{Algebraic
Geometry over Groups. 1. Algebraic Sets and Ideal Theory.} Journal of
Algebra. v.219, (1999), p. 16 -- 79.

[Ca] R. Camm. \textit{Simple free products.} J. London Math. Soc., \textbf{%
28,} 66-76, 1953.

[GSh] R.Gobel, S. Shelah. \textit{Radicals and Plotkin's problem concerning
geometrically equivalent groups.} Proc. Amer. Math. Soc., v. 130, (2002), p.
673 -- 674.

[GL] Gratzer G., Lakser H. \textit{A note on implicational class generated
by a class of structures,} Can. Math. Bull. (1974), v.16, n.4, p. 603 -- 605.

[Gv] Gvaramiya A. A. \textit{Quasi-varieties of many-sorted algebras,}
Theses of short reports in the international mathematical congress. Warsawa,
1983. Section 2, Algebra.

[Ma] Malcev A.I., \textit{Algebraic systems,} North Holland, 1973.

[MR] A.Myasnikov, V.Remeslennikov, \textit{Algebraic geometry over groups
II, Logical foundations,} J. of Algebra, \textbf{234:1} (2000) 225 -- 276.

[Pi] R. Pirs, \textit{Associative algebras,} Springer Verlag, 1982.

[Pl1] Plotkin B. \textit{Algebraic logic, varieties of algebras and
algebraic varieties,} Proc. Int. Alg. Conf., St. Petersburg, 1995,
St.Petersburg, 1999, p. 189 -- 271.

[Pl2] Plotkin B. \textit{Varieties of algebras and algebraic varieties.
Categories of algebraic varieties,} Siberian Advances in Mathematics,
v.7(2), (1997), p. 64 -- 97.

[Pl3] Plotkin B. \textit{Seven lectures on the Universal Algebraic Geometry,}
Preprint, {\ http:// arxiv:math, GM/0204245}, (2002), 87pp.

[Pl4] Plotkin B. \textit{Algebras with the same (algebraic) geometry,}
Proceedings of the International Conference on Mathematical Logic, Algebra
and Set Theory, dedicated to 100 anniversary of P.S.Novikov, Proceedings of
the Steklov Institute of Mathematics, MIAN, {\ v.242}, (2003), p. 17 -- 207.

[Pl5] Plotkin B. \textit{Action type logic and action type algebraic
geometry in the variety of group representations.} Manuscript.

[PPT] Plotkin B., Plotkin E., Tsurkov A. \textit{Geometrical equivalence of
groups,} Communications in Algebra. 27(8), 1999.

[PV] Plotkin B.I., Vovsi, S.M. \textit{Varieties of Group Representation},
Zinatne, Riga, 1983, (Russian).

[Vvs] Vovsi, S.M. \textit{Topics in Varieties of Group Representation},
Cambridge University Press, 1991.\pagebreak

\end{document}